\documentclass[journal]{IEEEtran/IEEEtran}
\ifCLASSINFOpdf
  \usepackage[pdftex]{graphicx}
\else
  \usepackage[dvips]{graphicx}
\fi

\usepackage[cmex10]{amsmath}
\usepackage{amsthm}

\usepackage{graphicx, bm, bbm, latexsym, amsfonts, epsfig, verbatim, epsfig, verbatim, enumerate, multirow, caption, graphicx, subcaption, algorithm, algpseudocode, hyperref, url, multicol, color, latexsym, amsfonts, epsfig, verbatim, tikz, cases, algcompatible, diagbox,multicol,fixmath,makecell,booktabs,comment}
\hypersetup{
	colorlinks=true, 
	breaklinks=true,
	urlcolor= black, 
	linkcolor= black, 
	citecolor=red, 
	pdftitle={},
	pdfauthor={},
}

\usepackage{graphicx, bbm, latexsym, amsfonts, epsfig, verbatim, amsmath, color, epsfig, verbatim, enumerate, multirow, caption, graphicx, subcaption, algorithm, algpseudocode}
\usepackage{url, multicol, color, latexsym, amsfonts, epsfig, verbatim, tikz, cases, bbm}
\usepackage[top=1.0in, bottom=1.0in, left=1.0in, right=1.0in]{geometry}
\usepackage{algcompatible}
\usepackage{diagbox,multicol,graphicx}
\usepackage{booktabs, caption, makecell}

\usepackage{threeparttable}
\usepackage{array}
\newcolumntype{L}[1]{>{\raggedright\let\newline\\\arraybackslash\hspace{0pt}}m{#1}}
\newcolumntype{C}[1]{>{\centering\let\newline\\\arraybackslash\hspace{0pt}}m{#1}}
\newcolumntype{R}[1]{>{\raggedleft\let\newline\\\arraybackslash\hspace{0pt}}m{#1}}

\DeclareMathOperator*{\argmin}{arg\,min}

\newtheorem{theorem}{Theorem}
\newtheorem{lemma}{Lemma}

\newtheorem{definition}{Definition}
\newtheorem{assumption}{Assumption}

\newtheorem{proposition}{Proposition}
\newtheorem{remark}{Remark}

\makeatletter
\newcommand{\subalign}[1]{%
  \vcenter{%
    \Let@ \restore@math@cr \default@tag
    \baselineskip\fontdimen10 \scriptfont\tw@
    \advance\baselineskip\fontdimen12 \scriptfont\tw@
    \lineskip\thr@@\fontdimen8 \scriptfont\thr@@
    \lineskiplimit\lineskip
    \ialign{\hfil$\m@th\scriptstyle##$&$\m@th\scriptstyle{}##$\hfil\crcr
      #1\crcr
    }%
  }%
}
\makeatother

\allowdisplaybreaks

\begin{document}
\title{A Privacy-Preserving Distributed Control of Optimal Power Flow
}

\author{Minseok~Ryu 
  and~Kibaek~Kim,~\IEEEmembership{Member,~IEEE,}
\thanks{M. Ryu and K. Kim are with the Mathematics and Computer Science Division, Argonne National Laboratory, Lemont, IL, USA (Contact: kimk@anl.gov).
This material is based upon work  supported by the U.S. Department of Energy, Office of Science, Advanced Scientific Computing Research, under Contract DE-AC02-06CH11357.
}
}

\maketitle

\begin{abstract}  
  We consider a distributed optimal power flow formulated as an optimization problem that maximizes a nondifferentiable concave function.
  Solving such a problem by the existing distributed algorithms can lead to data privacy issues because the solution information exchanged within the algorithms can be utilized by an adversary to infer the data.
  To preserve data privacy, in this paper we propose a differentially private projected subgradient (DP-PS) algorithm that includes a solution encryption step. 
  We show that a sequence generated by DP-PS converges in expectation, in probability, and with probability $1$.
  Moreover, we show that the rate of convergence in expectation is affected by a target privacy level of DP-PS chosen by the user.
  We conduct numerical experiments that demonstrate the convergence and data privacy preservation of DP-PS.
\end{abstract}


%
\IEEEpeerreviewmaketitle

\vspace{-4mm}
\section{Introduction}
Optimal power flow (OPF) is an important problem in reliably and economically operating electric grids. Currently, the problem is solved by
independent system operators  in a centralized manner.
Recently, however, distributed OPF has been spotlighted as a result of the introduction of microgrids with energy storage \cite{6469193} and increasing penetrations of distributed energy resources \cite{molzahn2017survey}.
Distributed OPF consists of (i) a set of OPF subproblems defined for each zone of the power grid and (ii) consensus constraints that link the subproblems. 
Distributed OPF can be solved by the existing distributed algorithms (e.g., \cite{sun2013fully,mhanna2018adaptive,mhanna2018component,sun2020two}), which do not require sharing private data information (e.g., demand data from each zone) but send the local solutions to the central machine.
Unfortunately, an adversary may be able to estimate the data based on the solutions (e.g., reverse engineering \cite{shokri2017membership}), thus  motivating the need for solution encryption.

Differential privacy (DP) is a randomization technique that guarantees the existence of multiple datasets with similar probabilities of resulting in the encrypted solution, thus preserving data privacy~\cite{dwork2006calibrating}.
A differentially private algorithm is an algorithm that incorporates  differential privacy for preserving data privacy during the algorithmic process \cite{dwork2014algorithmic}.
Several DP algorithms have been proposed to solve various distributed optimization problems.
For example, (i) a DP alternating direction method of multipliers (ADMM) was proposed for solving a distributed empirical risk minimization problem \cite{zhang2016dynamic, zhang2016dual, huang2019dp} and a distributed DC OPF \cite{dvorkin2019differentially}, and (ii) a DP stochastic gradient descent (SGD) method was proposed for solving a classification problem \cite{song2013stochastic}, a resource allocation problem \cite{han2016differentially}, and deep neural networks \cite{abadi2016deep}.

We define target privacy level (TPL) as a user parameter for  DP algorithms to control data privacy.
While guaranteeing stronger privacy, increasing TPL of the DP algorithms may affect the convergence.
For example, DP-ADMM with higher TPL is shown to find suboptimal solutions, implying the need for a trade-off between data privacy and solution quality~\cite{zhang2016dual, dvorkin2019differentially}.
On the other hand, the numerical results in~\cite{song2013stochastic} show that the solution accuracy from DP-SGD may be close to that of non-private SGD.
Also, Huang et al. \cite{huang2019dp} report  numerical experiments 
that DP-SGD has good noise resilience  compared with that of DP-ADMM, but it converges slowly.  
While numerical evidence has been demonstrated for the trade-off between the convergence of DP algorithms and TPL, only a few studies (e.g., DP-ADMM in~\cite{huang2019dp}) develop the theoretical links.

In this paper we present a DP projected subgradient (PS) algorithm for solving a distributed OPF while preserving data privacy, and we study how TPL affects the convergence of DP-PS theoretically and numerically.
We first formulate the distributed second-order conic (SOC) and alternating current (AC) OPF (see, e.g., \cite{coffrin2015qc, coffrin2016strengthening, kocuk2016strong}) based on dual decomposition by taking the Lagrangian relaxation with respect to the consensus constraints, where supergradients are computed by solving the OPF subproblems in parallel.
Moreover, in order to guarantee data privacy, the supergradients exchanged within the algorithm are systematically randomized by adding random noise extracted from a Laplace distribution.
Under three  rules of specifying search direction and step size,
we show that a sequence generated by DP-PS converges in expectation, in probability, and with probability $1$.
In particular, we show that the convergence complexity is affected by a constant factor only as TPL increases.
{\color{black} In summary, this paper answers the following questions:
\begin{itemize}
    \item How can we preserve the privacy of the data communicated in the distributed OPF setting? 
    \item Can DP guarantee data privacy?
    \item What are the implications of adding the DP technique to the distributed OPF setting?
    \item Can our technical development be numerically demonstrated?
\end{itemize}
} 
The remainder of the paper is organized as follows.
In Section \ref{sec:prelim}, we present a distributed OPF problem.
In Section \ref{sec:distributed_algo}, we describe a differentially private control with the proposed DP-PS.
In Section \ref{sec:algorithm}, we study the convergence of DP-PS.
We conduct case studies in Section \ref{sec:numerical} and summarize our conclusions in Section \ref{sec:conclusion}.
We denote by $\mathbb{N}$ a set of natural numbers.
For $A \in \mathbb{N}$, we define $[A]:=\{1,\ldots, A\}$.
We use $\langle \cdot, \cdot \rangle$ and $\| \cdot \|$ to denote the scalar product and the Euclidean norm.

\vspace{-4mm}
\section{Distributed Optimal Power Flow} \label{sec:prelim}

We depict a power network by a graph $(\mathcal{N},
\mathcal{L})$, where $\mathcal{N}$ is a set of buses and
$\mathcal{L}$ is a set of lines. For every line $\ell_{ij}
\in \mathcal{L}$, where $i$ is a \textit{from bus}  and $j$
is a \textit{to bus} of line $\ell$, we are given line
parameters, including bounds $[\underline{\theta}_{ij},
\overline{\theta}_{ij}]$ on voltage angle difference,
thermal limit $\overline{s}_{\ell}$, resistance $r_{\ell}$,
reactance $x_{\ell}$, impedance $z_{\ell} := r_{\ell} +
\textbf{i} x_{\ell}$, line charging susceptance
$b^{\text{\tiny c}}_{\ell}$, tap ratio $\tau_{\ell}$, phase
shift angle $\theta^{\text{\tiny s}}_{\ell}$, and admittance
matrix $Y_{\ell}$, namely,

\vspace{-3mm}
\begin{footnotesize}
\begin{align*}
Y_{\ell}
:= &
\begin{bmatrix}
Y^{{\text{\tiny ff}}}_{\ell} & Y^{{\text{\tiny ft}}}_{\ell} \\
Y^{{\text{\tiny tf}}}_{\ell} & Y^{{\text{\tiny tt}}}_{\ell}
\end{bmatrix}
=
\begin{bmatrix}
  (z^{-1}_{\ell} + \textbf{i} \frac{b^{\text{\tiny c}}_{\ell}}{2}) \frac{1}{\tau_{\ell}^2}
& -z^{-1}_{\ell} \frac{1}{\tau_{\ell} e^{- \textbf{i} \theta^{\text{\tiny s}}_{\ell}}}  \\
-z^{-1}_{\ell} \frac{1}{\tau_{\ell} e^{ \textbf{i} \theta^{\text{\tiny s}}_{\ell}}}
& z^{-1}_{\ell} + \textbf{i} \frac{b^{\text{\tiny c}}_{\ell}}{2}
\end{bmatrix}
,
\end{align*}
\end{footnotesize}
\noindent
$G^{\text{\tiny cf}}_{\ell} := \Re (Y^{\text{\tiny
ff}}_{\ell})$, $B^{\text{\tiny cf}}_{\ell} := \Im
(Y^{\text{\tiny ff}}_{\ell})$, $G^{\text{\tiny f}}_{\ell} :=
\Re (Y^{\text{\tiny ft}}_{\ell})$, $B^{\text{\tiny
f}}_{\ell} := \Im (Y^{\text{\tiny ft}}_{\ell})$,
$G^{\text{\tiny ct}}_{\ell} := \Re (Y^{\text{\tiny
tt}}_{\ell})$, $B^{\text{\tiny ct}}_{\ell} := \Im
(Y^{\text{\tiny tt}}_{\ell})$, $G^{\text{\tiny t}}_{\ell} =
\Re (Y^{\text{\tiny tf}}_{\ell})$, and $B^{\text{\tiny
t}}_{\ell} := \Im (Y^{\text{\tiny tf}}_{\ell})$. 
For every bus $i \in \mathcal{N}$, we are given bus parameters,
including bounds $[\underline{v}_i, \overline{v}_i]$ on
voltage magnitude, active (resp., reactive) power demand
$p^{\text{\tiny d}}_i$ (resp., $q^{\text{\tiny d}}_i$),
shunt conductance $g^{\text{\tiny s}}_i$, and shunt
susceptance $b^{\text{\tiny s}}_i$. Furthermore, for every
$i \in \mathcal{N}$, we define subsets
$\mathcal{L}^{\text{\tiny F}}_i := \{ \ell_{ij} : j \in
\mathcal{N}, \ell_{ij} \in \mathcal{L} \} $ and
$\mathcal{L}^{\text{\tiny T}}_i := \{ \ell_{ji} : j \in
\mathcal{N}, \ell_{ji} \in \mathcal{L} \}$ of $\mathcal{L}$
and a set of generators $\mathcal{G}_i$. For every generator
$g \in \mathcal{G}_i$, we are given generator parameters,
including bounds $[\underline{p}^{\text{\tiny G}}_g,
\overline{p}^{\text{\tiny G}}_g]$ (resp.,
$[\underline{q}^{\text{\tiny G}}_g,
\overline{q}^{\text{\tiny G}}_g]$) on the amounts of active
(resp., reactive) power generation and coefficients
($c_{1,g}$, $c_{2,g}$) of the quadratic generation cost
function.

Next we present decision variables. 
For every line $\ell_{ij} \in \mathcal{L}$, we
denote active (resp., reactive) power flow along line $\ell$
by $p^{\text{\tiny F}}_{\ell}$, $p^{\text{\tiny T}}_{\ell}$
(resp., $q^{\text{\tiny F}}_{\ell}$, $q^{\text{\tiny
T}}_{\ell}$).
For every $i \in \mathcal{N}$, we denote the complex voltage by
$V_i  = v^{\text{\tiny R}}_i + \textbf{i} v^{\text{\tiny
I}}_i$, and we introduce the following auxiliary variables:

\vspace{-3mm}
\begin{small}
\begin{align}
& w^{\text{\tiny RR}}_{ij} = v^{\text{\tiny R}}_i v^{\text{\tiny R}}_j, \ \
w^{\text{\tiny II}}_{ij} = v^{\text{\tiny I}}_i v^{\text{\tiny I}}_j, \ \
w^{\text{\tiny RI}}_{ij} = v^{\text{\tiny R}}_i v^{\text{\tiny I}}_j, \ \forall j \in \mathcal{N}. \label{AC_Linking}
\end{align}
\end{small}
\noindent
For every generator $g \in \mathcal{G}_i$, we
denote the amounts of active (resp., reactive) power
generation by $p^{\text{\tiny G}}_g$ (resp., $q^{\text{\tiny
G}}_g$).
In the following, we present a SOC OPF formulation:

\vspace{-3mm}
\begin{footnotesize}
\begin{subequations}
\label{ACOPF-rect}
\begin{align}
& \min \ \sum_{i \in \mathcal{N}} \sum_{g \in \mathcal{G}_i} \Big( c_{1,g} p^{\text{\tiny G}}_g + c_{2,g} (p^{\text{\tiny G}}_g )^2 \Big) \label{ACOPF-rect-0} \\
& \mbox{subject to} \nonumber \\ 
& \forall \ell_{ij} \in \mathcal{L}: \nonumber \\
& \ p^{\text{\tiny F}}_{\ell} = G^{\text{\tiny f}}_{\ell} (w^{\text{\tiny RR}}_{ij} + w^{\text{\tiny II}}_{ij}) +  B^{\text{\tiny f}}_{\ell} (w^{\text{\tiny RI}}_{ji} + w^{\text{\tiny RI}}_{ij}) + G^{\text{\tiny cf}}_{\ell} ( w^{\text{\tiny RR}}_{ii}+ w^{\text{\tiny II}}_{ii} ) , \label{ACOPF-rect-1} \\
& \ q^{\text{\tiny F}}_{\ell} = G^{\text{\tiny f}}_{\ell} (w^{\text{\tiny RI}}_{ji} + w^{\text{\tiny RI}}_{ij}) - B^{\text{\tiny f}}_{\ell} (w^{\text{\tiny RR}}_{ij} + w^{\text{\tiny II}}_{ij}) -B^{\text{\tiny cf}}_{\ell} ( w^{\text{\tiny RR}}_{ii}+ w^{\text{\tiny II}}_{ii} ),  \label{ACOPF-rect-2}\\
& \ p^{\text{\tiny T}}_{\ell} = G^{\text{\tiny t}}_{\ell} (w^{\text{\tiny RR}}_{ji} + w^{\text{\tiny II}}_{ji}) +  B^{\text{\tiny t}}_{\ell} (w^{\text{\tiny RI}}_{ij} + w^{\text{\tiny RI}}_{ji}) + G^{\text{\tiny ct}}_{\ell} ( w^{\text{\tiny RR}}_{jj}+ w^{\text{\tiny II}}_{jj} ), \label{ACOPF-rect-3} \\
& \ q^{\text{\tiny T}}_{\ell} = G^{\text{\tiny t}}_{\ell} (w^{\text{\tiny RI}}_{ij} + w^{\text{\tiny RI}}_{ji}) - B^{\text{\tiny t}}_{\ell} (w^{\text{\tiny RR}}_{ji} + w^{\text{\tiny II}}_{ji}) -B^{\text{\tiny ct}}_{\ell} ( w^{\text{\tiny RR}}_{jj}+ w^{\text{\tiny II}}_{jj} ), \label{ACOPF-rect-4} \\
& \ (p^{\text{\tiny F}}_{\ell})^2 + (q^{\text{\tiny F}}_{\ell})^2 \leq (\overline{s}_{\ell})^2, \ \ \ (p^{\text{\tiny T}}_{\ell})^2 + (q^{\text{\tiny T}}_{\ell})^2 \leq (\overline{s}_{\ell})^2, \label{ACOPF-rect-5} \\
& \ w^{\text{\tiny RI}}_{ji} - w^{\text{\tiny RI}}_{ij} \in [\tan(\underline{\theta}_{ij})(w^{\text{\tiny RR}}_{ij} + w^{\text{\tiny II}}_{ij}), \tan(\overline{\theta}_{ij})(w^{\text{\tiny RR}}_{ij} + w^{\text{\tiny II}}_{ij})], \label{ACOPF-rect-6} \\
& \forall i \in \mathcal{N}: \nonumber \\
& \ \sum_{\ell \in \mathcal{L}^{\text{\tiny F}}_i} p^{\text{\tiny F}}_{\ell} + \sum_{\ell \in \mathcal{L}^{\text{\tiny T}}_i}  p^{\text{\tiny T}}_{\ell} = \sum_{g \in \mathcal{G}_i} p^{\text{\tiny G}}_g - p^{\text{\tiny d}}_i - g^{\text{\tiny s}}_i (w^{\text{\tiny RR}}_{ii}+w^{\text{\tiny II}}_{ii}), \label{ACOPF-rect-7} \\
& \ \sum_{\ell \in \mathcal{L}^{\text{\tiny F}}_i} q^{\text{\tiny F}}_{\ell} + \sum_{\ell \in \mathcal{L}^{\text{\tiny T}}_i}  q^{\text{\tiny T}}_{\ell} = \sum_{g \in \mathcal{G}_i} q^{\text{\tiny G}}_g - q^{\text{\tiny d}}_i + b^{\text{\tiny s}}_i (w^{\text{\tiny RR}}_{ii}+w^{\text{\tiny II}}_{ii}), \label{ACOPF-rect-8} \\
& \ w^{\text{\tiny RR}}_{ii} + w^{\text{\tiny II}}_{ii}  \in [\underline{v}_i^2, \overline{v}_i^2], \label{ACOPF-rect-9} \\
& \forall i \in \mathcal{N}, \forall g \in \mathcal{G}_i: 
 p^{\text{\tiny G}}_g \in [\underline{p}^{\text{\tiny G}}_g, \overline{p}^{\text{\tiny G}}_g], \ \ \ q^{\text{\tiny G}}_g \in [\underline{q}^{\text{\tiny G}}_g, \overline{q}^{\text{\tiny G}}_g], \label{ACOPF-rect-10} \\
& {\color{black} \forall \ell_{ij} \in \mathcal{L}:} \nonumber \\
&  \ (w^{\text{\tiny RR}}_{ij}+w^{\text{\tiny II}}_{ij})^2 + (w^{\text{\tiny RI}}_{ji}-w^{\text{\tiny RI}}_{ij})^2 + \Big( \frac{w^{\text{\tiny RR}}_{ii}+w^{\text{\tiny II}}_{ii} - w^{\text{\tiny RR}}_{jj}-w^{\text{\tiny II}}_{jj} }{2} \Big)^2  \nonumber\\
& \ \leq \Big( \frac{w^{\text{\tiny RR}}_{ii}+w^{\text{\tiny II}}_{ii} + w^{\text{\tiny RR}}_{jj}+w^{\text{\tiny II}}_{jj} }{2} \Big)^2, \label{ACOPF-rect-11}
\end{align}
\end{subequations}
\end{footnotesize}
where \eqref{ACOPF-rect-0} is to minimize the generation cost, 
\eqref{ACOPF-rect-1}--\eqref{ACOPF-rect-4}
represent power flow, 
\eqref{ACOPF-rect-5} represent line thermal limit, 
\eqref{ACOPF-rect-6} represent bounds on voltage angle differences,
\eqref{ACOPF-rect-7}--\eqref{ACOPF-rect-8} represent power
balance, 
\eqref{ACOPF-rect-9} represent bounds on voltage magnitudes,
\eqref{ACOPF-rect-10} represent bounds on power generation, and 
\eqref{ACOPF-rect-11} represent SOC constraints that ensure linking between auxiliary variables.

{\color{black}
We remark that this paper uses the SOC OPF formulation for example. The technical development and results should remain true with any convex relaxation of the OPF problem.
For example, one can introduce SOCP strengthening techniques \cite{kocuk2016strong, bynum2018strengthened}, semidefinite programming relaxation, or quadratic convex relaxation \cite{coffrin2015qc}. 
In this work, however, we focus on solving one of the convex relaxation techniques, SOC OPF \eqref{ACOPF-rect} in a distributed and privacy-preserving manner.
}

We decompose the network into several zones indexed by $\mathcal{Z}:= \{1, \ldots,
Z\}$. Specifically, we split a set $\mathcal{N}$ of buses
into subsets $\{ \mathcal{N}_z\}_{z \in \mathcal{Z}}$ such
that $\mathcal{N} = \cup_{z \in \mathcal{Z}} \mathcal{N}_z$
and $\mathcal{N}_z \cap \mathcal{N}_{z'} = \emptyset$ for
$z,z' \in \mathcal{Z}: z \neq z'$. For each zone $z \in
\mathcal{Z}$ we define a line set $\mathcal{L}_z:= \cup_{i
\in \mathcal{N}_z} \big( \mathcal{L}^{\text{\tiny F}}_i \cup
\mathcal{L}^{\text{\tiny T}}_i \big)$; an extended node set
$\mathcal{V}_z := \cup_{i \in \mathcal{N}_z} \mathcal{A}_i$,
where $\mathcal{A}_i$ is a set of adjacent buses of $i$; and
a set of cuts $\mathcal{C}_z = \cup_{z' \in
\mathcal{Z}\setminus \{z\} } (\mathcal{L}_z \cap
\mathcal{L}_{z'})$. Note that $\{ \mathcal{N}_z \}_{z \in
\mathcal{Z}}$ is a collection of disjoint sets, while $\{
\mathcal{L}_z \}_{z \in \mathcal{Z}}$ and $\{ \mathcal{V}_z
\}_{z \in \mathcal{Z}}$ are not. 
Using these notations,
we rewrite problem \eqref{ACOPF-rect}  as

\vspace{-3mm}
\begin{small}
\begin{subequations}
\label{model:ACOPF_matrix}
\begin{align}
\min \ & \sum_{z \in \mathcal{Z}} f_z(x_z) \\
\mbox{s.t.} \ & (x_z, y_z) \in \mathcal{F}_z (\bar{D}_z), \ \forall z \in \mathcal{Z}, \label{model:ACOPF_matrix-1} \\
& \phi_i = y_{zi}, \ \forall z \in \mathcal{Z}, \forall i \in C(z), \label{model:ACOPF_matrix-2} \\
& \phi_i \in \mathbb{R}, \ \forall i \in \mathcal{C}, \label{model:ACOPF_matrix-3}
\end{align}
\end{subequations}
\end{small}
where

\vspace{-3mm}
\begin{small}
\begin{align*}
  & x_z \leftarrow \big\{p^{\text{\tiny F}}_{z\ell}, q^{\text{\tiny F}}_{z\ell}, p^{\text{\tiny T}}_{z\ell}, q^{\text{\tiny T}}_{z\ell},
  w^{\text{\tiny RR}}_{zij}, w^{\text{\tiny II}}_{zij}, w^{\text{\tiny RI}}_{zij}, w^{\text{\tiny RI}}_{zji} \big\}_{\ell_{ij} \in \mathcal{L}_z \setminus \mathcal{C}_z} \\
  & \hspace{8mm} \cup \big\{ v^{\text{\tiny R}}_{zi}, v^{\text{\tiny I}}_{zi} \big\}_{i \in \mathcal{V}_z} \cup \big\{ p^{\text{\tiny G}}_g, q^{\text{\tiny G}}_g \big\}_{i \in \mathcal{N}_z, g \in \mathcal{G}_i}, \\  
& y_z \leftarrow \big\{p^{\text{\tiny F}}_{z\ell}, q^{\text{\tiny F}}_{z\ell}, p^{\text{\tiny T}}_{z\ell}, q^{\text{\tiny T}}_{z\ell},
w^{\text{\tiny RR}}_{zij}, w^{\text{\tiny II}}_{zij}, w^{\text{\tiny RI}}_{zij}, w^{\text{\tiny RI}}_{zji} \big\}_{\ell_{ij} \in \mathcal{C}_z}, \\
& \phi \leftarrow \cup_{z \in \mathcal{Z}} \big\{ p^{\text{\tiny F}}_{\ell}, q^{\text{\tiny F}}_{\ell}, p^{\text{\tiny T}}_{\ell}, q^{\text{\tiny T}}_{\ell},
w^{\text{\tiny RR}}_{ij}, w^{\text{\tiny II}}_{ij}, w^{\text{\tiny RI}}_{ij}, w^{\text{\tiny RI}}_{ji} \big\}_{\ell_{ij} \in \mathcal{C}_z}, 
\end{align*}
\end{small}
$C(z)$ is an index set that indicates each element of $y_z$,
$\mathcal{C} := \cup_{z \in \mathcal{Z}} C(z)$ is an index set of consensus variable {\color{black} $\phi$}, $f_z(x_z):= \sum_{i \in \mathcal{N}_z} \sum_{g \in \mathcal{G}_i} \big( c_{1,g} p^{\text{\tiny G}}_g + c_{2,g} (p^{\text{\tiny G}}_g )^2 \big)$, $\bar{D}_z := \{ p^{\text{\tiny d}}_l
\}_{l \in \mathcal{N}_z}$ is a given demand vector, $\mathcal{F}_z(\bar{D}_z) := \{(x_z, y_z): \eqref{ACOPF-rect-1}-\eqref{ACOPF-rect-6}, \forall \ell_{ij} \in \mathcal{L}_z; \ \eqref{ACOPF-rect-7}, \eqref{ACOPF-rect-8}, \forall i \in \mathcal{N}_z; \ \eqref{ACOPF-rect-9},\forall i \in \mathcal{V}_z; \ \eqref{ACOPF-rect-10},\forall i \in \mathcal{N}_z,\forall g \in \mathcal{G}_i; \ \eqref{ACOPF-rect-11}, \forall i \in \mathcal{V}_z,\forall j\in \mathcal{V}_z \} $ is a convex feasible region defined for each zone, and \eqref{model:ACOPF_matrix-2} represents the consensus constraints:

\vspace{-3mm}
\begin{small}
\begin{align*}
& \forall z \in \mathcal{Z}, \ell_{ij} \in \mathcal{C}_z: \nonumber \\
& \ p^{\text{\tiny F}}_{\ell} = p^{\text{\tiny F}}_{z \ell}, \ \ p^{\text{\tiny T}}_{\ell} = p^{\text{\tiny T}}_{z \ell}, \ \
q^{\text{\tiny F}}_{\ell} = q^{\text{\tiny F}}_{z \ell}, \ \ q^{\text{\tiny T}}_{\ell} = q^{\text{\tiny T}}_{z \ell}, \\
& \ w^{\text{\tiny RR}}_{ij} = w^{\text{\tiny RR}}_{zij}, \  w^{\text{\tiny II}}_{ij} = w^{\text{\tiny II}}_{zij}, \ 
w^{\text{\tiny RI}}_{ij} = w^{\text{\tiny RI}}_{zij}, \  w^{\text{\tiny RI}}_{ji} = w^{\text{\tiny RI}}_{zji}.  
\end{align*}
\end{small} 
Note that the consensus constraints with respect to $w^{\text{\tiny RR}}_{ij},w^{\text{\tiny II}}_{ij},w^{\text{\tiny RI}}_{ij},w^{\text{\tiny RI}}_{ji}$ are redundant but numerically beneficial.
By introducing a dual vector $\lambda:=\{\lambda_{zi}\}_{z \in \mathcal{Z}, i \in C(z)}$
associated with constraints \eqref{model:ACOPF_matrix-2}, one can construct a Lagrangian dual problem:

\begin{subequations}
\label{model:Lagrangian_ACOPF_matrix}
\vspace{-3mm}
\begin{small}
\begin{align}
\max_{\lambda \in \Lambda} \ & \Big\{ H(\lambda) := \sum_{z \in \mathcal{Z}} h_z(\lambda_z) \Big\},
\end{align}
\end{small}
where $\Lambda := \big\{ \lambda: \ \sum_{z \in F(i)} \lambda_{zi} = 0, \forall i \in \mathcal{C} \big\}$, $F(i) := \{z \in \mathcal{Z}: i \in C(z) \}$ is a set of zones for every  $i \in \mathcal{C}$ and $h_z(\lambda_z)$ is the optimal value of the subproblem:

\vspace{-3mm}
\begin{small}
\begin{align}
& \min_{ (x_z,y_z) \in \mathcal{F}_z(\bar{D}_z)} \ f_z(x_z) + \sum_{i \in C(z)} \lambda_{zi} y_{zi}. \label{model:nonconvex_subproblem}
\end{align}
\end{small}
\end{subequations}
{\color{black}
We emphasize that, for a given $\lambda$, evaluating $H(\lambda)$ can be done by solving the subproblem \eqref{model:nonconvex_subproblem} in parallel.
}
Let $\lambda^*$ be a maximizer of the nondifferentiable concave function $H(\lambda)$.
Then $H(\lambda^*)$ is the optimal value of \eqref{model:ACOPF_matrix} by the strong duality from the convexity of $\mathcal{F}_z(\bar{D}_z)$.
\begin{remark} \label{remark:ACOPF_model}
  If \eqref{ACOPF-rect-11} are replaced with \eqref{AC_Linking}, then \eqref{ACOPF-rect} is a rectangular formulation of AC OPF. 
  In this case, \eqref{model:Lagrangian_ACOPF_matrix} may not provide a solution that satisfies \eqref{model:ACOPF_matrix-2} because of the nonconvexity of $\mathcal{F}_z(\bar{D}_z)$.
\end{remark}
\begin{remark} \label{remark:bound_y}
There exists $y^\text{\tiny L}_i, y^\text{\tiny U}_i \in \mathbb{R}$ such that $y_{zi} \in [y^\text{\tiny L}_i, y^\text{\tiny U}_i], \ \forall i \in \mathcal{C}, \forall z \in F(i)$.
\end{remark}

\vspace{-6mm}
\section{Differentially Private Control} \label{sec:distributed_algo}
The Lagrangian dual problem \eqref{model:Lagrangian_ACOPF_matrix} can be solved by any nonsmooth convex optimization algorithms.
In this paper we consider the PS algorithm,

\vspace{-3mm}
\begin{small}
\begin{align}
  \lambda^{k+1} = \text{Proj}_{\Lambda}( \lambda^k + \alpha_k y^k), \ \forall k \in [K], \label{subgradient_projection_algorithm}
\end{align}
\end{small}
where $\text{Proj}_{\Lambda}(\cdot)$ represents the orthogonal projection onto $\Lambda$, $\alpha_k$ is a step size, $y^k$ is a search direction, and $K$ is the total number of iterations.

\subsubsection{Motivating Example (Data Leakage)} \label{sec:motivating}

{\color{black}
Throughout the paper, we consider a {\color{black} hypothetically strong adversary that can access every but} private load data of a control zone in a power system and tries to infer the data by intercepting the communication among the control zones.
In this example, we demonstrate that the existing distributed algorithms are susceptible to inference attacks (e.g., \cite{shokri2017membership}).
Specifically, the adversary can intercept the communication data $\{y^k_z\}_{k=1}^{K}$ of the PS algorithm \eqref{subgradient_projection_algorithm} and try to infer private demand data $\bar{D}_z$ in \eqref{model:nonconvex_subproblem}.
Such inference attack can be easily conducted by solving an adversary problem as in \cite{dvorkin2019differentially}.
The adversary problem is described as follows.

Let $\mathcal{K}$ be a set of PS iterations observed by an adversary who aims to infer a demand data at node $\hat{l}$ in zone $\hat{z}$, namely $\bar{D}_{\hat{z} \hat{l}}$.
We {\color{black} assume that the strong} adversary knows
(i) all the demand information except $\bar{D}_{\hat{z}\hat{l}}$,
(ii) all the topological information of zone $\hat{z}$, and
(iii) the exchanged supergradient $\{\hat{y}^k_{\hat{z}} \}_{k \in \mathcal{K}}$ and the local solution $\{ \hat{x}^k_{\hat{z}} \}_{k \in \mathcal{K}}$.
{\color{black}
Our assumption is justified as to give the most advantages to the adversary, which can be considered as the worst-case data leakage scenario to the pravacy-preserving control.
}



\vspace{-3mm}
\begin{footnotesize}
\begin{align}
& \min_{\substack{D_{\hat{z}\hat{l}}, x^k_{\hat{z}}, y^k_{\hat{z}}}} \ \sum_{k \in \mathcal{K}} f_{\hat{z}}(x^k_{\hat{z}}) + \Gamma \big\{ \| x^k_{\hat{z}} - \hat{x}^k_{\hat{z}} \|^2 + \| y^k_{\hat{z}} - \hat{y}^k_{\hat{z}} \|^2 \big\}  \label{Adversarial} \\
& \mbox{subject to} \ \  \forall k \in \mathcal{K}: \nonumber \\
& \ \eqref{ACOPF-rect-1}-\eqref{ACOPF-rect-6}, \forall \ell_{ij} \in \mathcal{L}_{\hat{z}}; \ \ \eqref{ACOPF-rect-8},\forall i \in \mathcal{N}_{\hat{z}}; \ \eqref{ACOPF-rect-9}, \forall i \in \mathcal{V}_{\hat{z}}; \nonumber \\
&  \ \eqref{ACOPF-rect-10}, \forall i \in \mathcal{N}_{\hat{z}}, \forall g \in \mathcal{G}_i; \ \eqref{ACOPF-rect-11}, \forall i \in \mathcal{V}_{\hat{z}}, \forall j \in \mathcal{V}_{\hat{z}}; \nonumber\\
& \ \sum_{\ell \in \mathcal{L}^{\text{\tiny F}}_{\hat{l}}} p^{\text{\tiny F}}_{k \hat{z} \ell} + \sum_{\ell \in \mathcal{L}^{\text{\tiny T}}_{\hat{l}}}  p^{\text{\tiny T}}_{k \hat{z} \ell} = \sum_{g \in \mathcal{G}_{\hat{l}}} p^{\text{\tiny G}}_{kg} - D_{\hat{z}\hat{l}} - g^{\text{\tiny s}}_{\hat{l}} (w^{\text{\tiny RR}}_{k\hat{z} \hat{l} \hat{l}}+w^{\text{\tiny II}}_{k\hat{z} \hat{l} \hat{l}}); \nonumber \\
& \ \sum_{\ell \in \mathcal{L}^{\text{\tiny F}}_{\hat{l}}} p^{\text{\tiny F}}_{k \hat{z} \ell} + \sum_{\ell \in \mathcal{L}^{\text{\tiny T}}_{\hat{l}}}  p^{\text{\tiny T}}_{k \hat{z} \ell} =  \sum_{g \in \mathcal{G}_{\hat{l}}} p^{\text{\tiny G}}_{kg} - p^{\text{\tiny d}}_{i} - g^{\text{\tiny s}}_{\hat{l}} (w^{\text{\tiny RR}}_{k\hat{z} \hat{l} \hat{l}}+w^{\text{\tiny II}}_{k\hat{z} \hat{l} \hat{l}}), \nonumber \\
& \ \ \ \ \forall i \in \mathcal{N}_{\hat{z}} \setminus \{\hat{l}\}, \nonumber
\end{align} 
\end{footnotesize}
where $\Gamma>0$ is a penalty parameter.
Note that $D_{\hat{z} \hat{l}}$ is a decision variable as well as $\{x^k_{\hat{z}}, y^k_{\hat{z}}\}_{k \in \mathcal{K}}$.
By solving \eqref{Adversarial}, the adversary aims to obtain the unknown demand $D_{\hat{z} \hat{l}}$ that minimizes the distance between $\{x^k_{\hat{z}}, y^k_{\hat{z}}\}_{k \in \mathcal{K}}$ and the solutions $\{\hat{x}^k_{\hat{z}}, \hat{y}^k_{\hat{z}}\}_{k \in \mathcal{K}}$ obtained from the PS algorithm.
With a sufficiently large $\Gamma$, the adversary can find the demand data at node $\hat{l}$ that produces $\{\hat{x}^k_{\hat{z}}, \hat{y}^k_{\hat{z}}\}_{k \in \mathcal{K}}$, thus identifying $\bar{D}_{\hat{z} \hat{l}}$.
As the cardinality of $\mathcal{K}$ increases, moreover, the accuracy of the demand estimated by \eqref{Adversarial} increases while sacrificing computation.
We denote by $\widehat{\mathcal{K}}$ a collection of various $\mathcal{K}$ and by ${D}_{\hat{z}\hat{l}}(\mathcal{K})$ a demand estimated by \eqref{Adversarial} with $\mathcal{K} \in \widehat{\mathcal{K}}$.
}

\begin{figure}[!h]
  \centering  
  \includegraphics[scale=0.2]{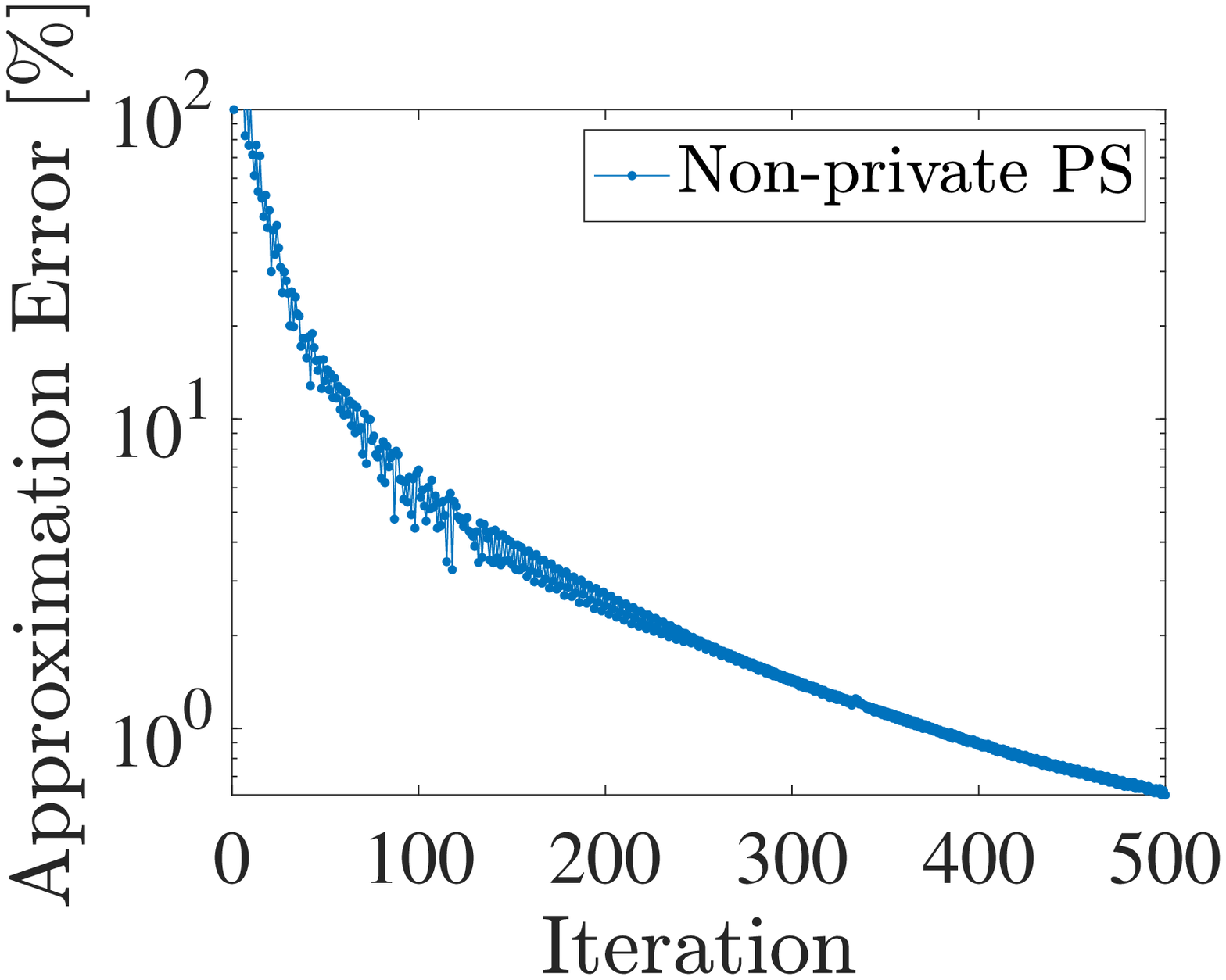}
  \includegraphics[scale=0.2]{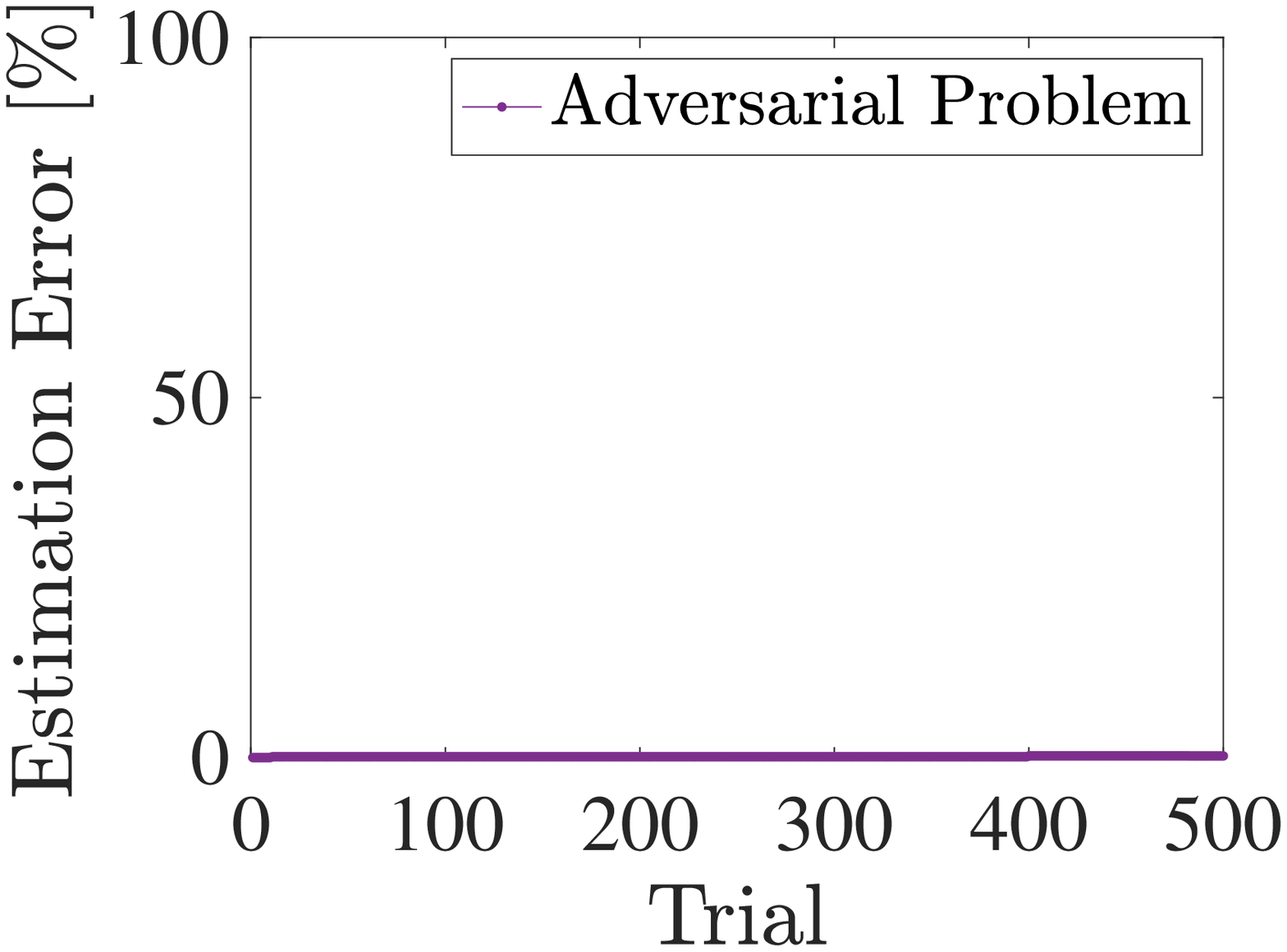}    
  \caption{Approximation error (left) of PS and  demand estimation error (right) of the adversarial problem.}
  \label{fig:case14_no_perturbation}
\end{figure}

We demonstrate the effectiveness of adversary problem~\eqref{Adversarial} by using an instance ``case 14'' from Matpower \cite{zimmerman2010matpower} with the decomposition into $3$ zones (see Table \ref{tab:decomposition}).
We solve the distributed OPF of~\eqref{model:Lagrangian_ACOPF_matrix} by using PS.
At each iteration $k$, an approximation error is measured as below and reported in Figure \ref{fig:case14_no_perturbation}:

\vspace{-3mm}
\begin{small}
\begin{align}
& \text{AE}_k = 100|Z^*-Z^k|/Z^*, \ \forall k \in [K], \label{AE} 
\end{align}
\end{small}
\noindent
where $Z^*$ is the optimal objective value and $Z^k$ are the objective values computed at the $k$th iteration of PS, respectively.

We consider the adversary who aims to estimate the demand at node $\hat{l} = 4$ in zone $\hat{z}=1$, namely, $\bar{D}_{14} = 47.8$ MW.
For every trial $\mathcal{K} \in \widehat{\mathcal{K}}$, we solve \eqref{Adversarial} and report in Figure \ref{fig:case14_no_perturbation} a demand estimation error:

\vspace{-3mm}
\begin{small}
\begin{align}
& \text{DE} (\mathcal{K}) :=  100 | \bar{D}_{\hat{z} \hat{l}} - {D}_{\hat{z}\hat{l}}(\mathcal{K}) |/\bar{D}_{\hat{z} \hat{l}}, \ \forall \mathcal{K} \in \widehat{\mathcal{K}}, \label{DE}
\end{align}
\end{small}
where $\widehat{K} \leftarrow \big\{ \{1\}, \ldots, \{K\} \big\}$ (various $\widehat{K}$ will be discussed in Section \ref{sec:numerical}).
Figure \ref{fig:case14_no_perturbation} shows that the adversary is highly likely to estimate $\bar{D}_{14}$, and hence this situation motivates the need for the solution encryption to preserve data privacy.

\subsubsection{Differential Privacy in PS}

The motivating example suggests that PS for solving \eqref{model:Lagrangian_ACOPF_matrix} might be vulnerable to data leakage.
To preserve data privacy, we introduce differential privacy (see \cite{dwork2014algorithmic} for more details).

\begin{definition}{\text{\bf ($\bar{\epsilon}$-differential privacy)}} \label{def:epsilon_differential_privacy}
A randomized function $\mathcal{R}$ that maps data $D$ to some random numbers gives $\bar{\epsilon}$-differential privacy if

\vspace{-3mm}
\begin{footnotesize}
\begin{align*}
\bigg\vert \ln \bigg( \frac{\mathbb{P} \{ \mathcal{R}(D') \in S \}}{\mathbb{P} \{ \mathcal{R}(D'') \in S \}} \bigg) \bigg\vert \leq \bar{\epsilon}, \ \forall (D',D'') \in \mathcal{D}_{\beta}, \forall S \subseteq \text{Range}(\mathcal{R}),
\end{align*}
\end{footnotesize}
\noindent
where $\bar{\epsilon} > 0$, the probability taken is over the coin tosses of $\mathcal{R}$, and $\mathcal{D}_{\beta}$ is a collection of two datasets $(D', D'')$ differing in one element by $\beta \in
\mathbb{R}_+$.
\end{definition}
\noindent
For small $\bar{\epsilon} \approx
ln(1+\bar{\epsilon})$, we have 
$\frac{\mathbb{P} \{ \mathcal{R}(D')
\in S \}}{\mathbb{P} \{ \mathcal{R}(D'') \in S \}} \in [1- \bar{\epsilon},
1+\bar{\epsilon}]$, which implies that distinguishing $D'$ from $D''$ based on $S$ becomes more difficult as $\bar{\epsilon}$ decreases.
To construct $\mathcal{R}(D)$ that ensures $\bar{\epsilon}$-differential privacy on data $D$, one can utilize a Laplace mechanism \cite{dwork2006calibrating}. More specifically, a query function $\mathcal{Q}:D \mapsto \mathbb{R}$ mapping data to \emph{true answer} is perturbed by adding Laplacian noise described in Definition \ref{def:laplace}.

\begin{definition}{\text{\bf (Laplacian noise)}} \label{def:laplace}
 Laplacian noise $\tilde{\xi} \in \mathbb{R}$ is a random variable following the Laplace distribution whose probability density function is
$L(\tilde{\xi}|  b) = \frac{1}{2  b} \exp \big( - \frac{\vert \tilde{\xi} \vert}{ b} \big)$ for $b>0$.
\end{definition}
\noindent
The randomized function $\mathcal{R}(D):=\mathcal{Q}(D) +\tilde{\xi}$ provides $\bar{\epsilon}$-differential privacy if $\tilde{\xi}$ is drawn from the Laplace distribution with $b = \max_{(D',D'') \in \mathcal{D}_{\beta}} |\mathcal{Q}(D') - \mathcal{Q}(D'')|/ \bar{\epsilon}$.

The main idea of DP-PS is to perturb $y^k$ with the noise $\tilde{\xi}^k$ such that

\vspace{-3mm}
\begin{small}
\begin{align}
\tilde{y}^k_{zi} \leftarrow  y^k_{zi} + \tilde{\xi}^k_{zi}, \ \forall z \in \mathcal{Z}, \forall i \in C(z), \label{Perturb_y}
\end{align}
\end{small}
\noindent
for every iteration $k$ of PS.
We describe the algorithmic steps of DP-PS in Algorithm \ref{algo:DPPSA}.
In line 3, we find a supergradient $y^k$ of the concave function $H$ at $\lambda^k$.
In lines 6--8, we generate the Laplacian noise $\tilde{\xi}^k$ and the noisy supergradient $\tilde{y}^k$.
In line 9, we update dual variables based on the step size $\alpha_k$ and the search direction $s^k(\tilde{y}^k)$ determined in advance (see Section \ref{sec:algorithm}).

\begin{algorithm}[!h]
  \caption{DP Projected Subgradient Algorithm}
  \begin{algorithmic}[1]
  \STATE Set $k \leftarrow 1$ and $\lambda^1 \leftarrow 0$.
  \FOR{$k \in \{1, \ldots, K\}$}   
  \STATE
  Given $\lambda^k$, find $y^k$ by solving \eqref{model:nonconvex_subproblem} in parallel.
  \STATE Store $H_{\text{best}} (\lambda^k) \leftarrow \max_{t \in [k]}
  \{ H( \lambda^t) \}$. 
  \STATE \textbf{\# Perturbation of $y^k$}
  {\color{black}
  \STATE Solve \eqref{model:Sensitivity_problem} to find $\{\bar{\Delta}^k_{zi} \}_{z \in \mathcal{Z}, i \in C(z)}$.
  \STATE Extract $\tilde{\xi}^k_{zi}$ from
  $L(\tilde{\xi}^k_{zi} \vert \bar{\Delta}^k_{zi} / \bar{\epsilon})$ in Definition \ref{def:laplace}.
  }
  \STATE Compute $\tilde{y}^k$ by \eqref{Perturb_y}.
  \STATE \textbf{\# Update dual variables}  

  \vspace{-3mm}
  \begin{small}
  \begin{align*}
    \lambda^{k+1}  \leftarrow \text{Proj}_{\Lambda} \big(\lambda^k + \alpha_k s^k(\tilde{y}^k) \big).
  \end{align*} 
  \end{small} 
  \ENDFOR
  \end{algorithmic}
  \label{algo:DPPSA}
\end{algorithm}

Now we describe how to generate the noise $\tilde {\xi}^k_{zi}$ in \eqref{Perturb_y} so that the $\bar{\epsilon}$-differential privacy in Definition \ref{def:epsilon_differential_privacy} on $\bar{D}$ is ensured.
First, we define a query function as follows:

\vspace{-3mm}
\begin{small}
\begin{align}
  \mathcal{Q}_{zi}^k: D_z \mapsto y^k_{zi}, \ \forall z \in \mathcal{Z}, \forall i \in C(z), \label{query_fn}
\end{align}
\end{small}
\noindent
where $y^k_{zi}$ is obtained by solving \eqref{model:nonconvex_subproblem} for given $\lambda^k$ and $D_z \in \mathbb{R}^{|\mathcal{N}_z|}$.
{\color{black}
Second, we draw $\tilde{\xi}^k_{zi}$ in \eqref{Perturb_y} from the Laplace distribution in Definition \ref{def:laplace} with $b=\bar{\Delta}^k_{zi} / \bar{\epsilon}$ and

\vspace{-3mm}
\begin{small}
\begin{align}
& \bar{\Delta}^k_{zi} := \max_{D'_z \in \widehat{\mathcal{D}}_{{\beta}}(\bar{D}_z) } | \mathcal{Q}^k_{zi} (D'_z) - \mathcal{Q}^k_{zi} (\bar{D}_z)|, \label{model:Sensitivity_problem}
\end{align}
\end{small}
where $\bar{D}_z$ is a given demand and $\widehat{\mathcal{D}}_{{\beta}}(\bar{D}_z)$ is a collection of $D'_z$ differing in one element from $\bar{D}_z$ by $\beta$, namely

\vspace{-3mm} 
\begin{small}
  \begin{align*}
  & \widehat{\mathcal{D}}_{\beta} (\bar{D}_z) := \cup_{l \in \mathcal{N}_z } \Big\{ D'_z  \in \mathbb{R}^{|\mathcal{N}_z|} : D'_{zj} = \bar{D}_{zj}, \ \forall j \in \mathcal{N}_z \setminus \{l\}, \\
  &   D'_{zl} \in [ \bar{D}_{zl} ( 1 - \beta), \ \bar{D}_{zl} ( 1+ \beta)] \Big\}.
  \end{align*}
\end{small}

\begin{theorem} \label{thm:privacy_guarantee}
  In the $k$th iteration of DP-PS, we consider
  
  \vspace{-3mm}
  \begin{small}
  \begin{align}
    \mathcal{R}^k_{zi}(D_z) := \mathcal{Q}^k_{zi} (D_z)+ \tilde{\xi}^k_{zi}, \ \forall z \in \mathcal{Z}, \forall i \in C(z), \label{Randomized_fn_ACOPF}
  \end{align}
  \end{small}
  \noindent
  where $\mathcal{Q}^k_{zi}$ is defined in \eqref{query_fn}, $\tilde{\xi}^k_{zi}$ is extracted from $L(\tilde{\xi}^k_{zi} \vert \bar{\Delta}^k_{zi} / \bar{\epsilon})$ in Definition \ref{def:laplace}, and $\bar{\Delta}^k_{zi}$ is from \eqref{model:Sensitivity_problem}.
  For all $z \in \mathcal{Z}$, $i \in C(z)$, we have
  
  \vspace{-3mm}
  \begin{small}    
    \begin{align}
      \bigg\vert \ln \bigg( \frac{\mathbb{P} \{ \mathcal{R}^k_{zi}(D'_z) \in S_k \}}{\mathbb{P} \{ \mathcal{R}^k_{zi}(\bar{D}_z) \in S_k  \}} \bigg) \bigg\vert \leq \bar{\epsilon}, \ \
      & \subalign{&\forall D'_z \in \widehat{\mathcal{D}}_{\beta}(\bar{D}_z), \\ &\forall S_k \subseteq \text{Range}(\mathcal{R}^k_{zi}),} \label{epsilon_dp}
    \end{align}  
  \end{small}
  \noindent
  where $\mathcal{R}^k_{zi}(\bar{D}_z)$ is equal to $\tilde{y}^k_{zi}$ in \eqref{Perturb_y}. This implies that $\bar{\epsilon}$-differential privacy on $\bar{D}$ is ensured in the $k$th iteration of DP-PS.
  Moreover, DP-PS with $K$ total iterations provides $\bar{\epsilon}$-differential privacy against an adversary observing the information exchanged during the entire process of DP-PS, if $\tilde{\xi}^k_{zi}$ is extracted from $L(\tilde{\xi}^k_{zi} \vert K  \bar{\Delta}^k_{zi} / \bar{\epsilon})$. 
  \end{theorem}
  \proof See Appendix \ref{apx-thm:privacy_guarantee}.
  \qed

  \begin{remark}{\text{\bf (Target Privacy Level)}} \label{remark:privacy_level}
    We denote by $1/\bar{\epsilon}$ a target privacy level (TPL) inproportional to the privacy loss $\bar{\epsilon}$.
    As TPL increases, the variance of $\tilde{\xi}^k_{zi}$, namely, $2(\bar{\Delta}^k_{zi} / \bar{\epsilon})^2$, increases for fixed $\bar{\Delta}^k_{zi}$ and stronger data privacy is achieved according to Definition \ref{def:epsilon_differential_privacy}.
  \end{remark}
}

\vspace{-3mm}
\section{Convergence of DP-PS}
\label{sec:algorithm}
In this section we study how TPL affects the convergence of DP-PS.
In Table \ref{tab:three_rules} we describe three rules for determining step size and search direction.
Note that the step size in Rule 1 is deterministic and \textit{square-summable but not summable} and that the step sizes in Rule 2 and Rule 3 are stochastic and affected by $\tilde{y}^k$, which are a variant of Polyak \cite{polyakintroduction} and CFM \cite{camerini1975improving}, respectively.
{\color{black} 
As compared with the previous work in \cite{huang2019dp} that proposes DP-ADMM guaranteeing the convergence in expectation, the DP-PS with the three rules provides the convergence in expectation and probability, and the almost sure convergence, which is a stronger result.
}

\begin{table}[htpb]
\centering
\caption{Three rules for DP-PS.}
\resizebox{\columnwidth}{!}{%
\begin{tabular}{|c|l|l|}
\hline  
Rule &  Step size  & Search direction \\ \hline
1  & \makecell[l]{$\alpha_k = a/k$, \\ where $ a > 0$.} & $s^k(\tilde{y}^k) := \tilde{y}^k$  \\ \hline
2  & $\alpha_k := \frac{{H}(\lambda^{\star}) - {H}(\lambda^{k})}{\| s^k(\tilde{y}^k)  \|^2}$ & $s^k(\tilde{y}^k) := \tilde{y}^k$ \\ \hline
3  & $\alpha_k := \frac{{H}(\lambda^{\star}) - {H}(\lambda^{k})}{\| s^k(\tilde{y}^k)  \|^2}$ & \makecell[l]{$s^k (\tilde{y}^k) := \tilde{y}^k + \zeta_k s^{k-1}(\tilde{y}^{k-1})$,  \\ $\zeta_k := \max \{ 0, - \chi_k \frac{\langle s^{k-1}(\tilde{y}^{k-1}), \tilde{y}^k \rangle}{\|s^{k-1}(\tilde{y}^{k-1})\|^2}  \}$, \\ where $s^0=0$ and $\chi_k \in [0,2]$}  \\ \hline
\end{tabular}
}
\label{tab:three_rules}
\end{table}

\begin{assumption} \label{assump:dual_bound}
$\Lambda$ is compact and $\lambda^* \in \Lambda$ maximizes $H$.
\end{assumption}

{\color{black} 
\begin{remark} \label{remark:noise}
For the Laplacian noise $\tilde{\xi}^k_{zi}$ in \eqref{Perturb_y}, we notice that there exists 
$\tilde{\xi}^{\text{\tiny U}}_{zi}( \bar{\epsilon}) \in \mathbb{R}_+$, which increases as $\bar{\epsilon}$ decreases, such that $\tilde{\xi}^k_{zi} \in [-\tilde{\xi}^{\text{\tiny U}}_{zi}( \bar{\epsilon} ), \tilde{\xi}^{\text{\tiny U}}_{zi}( \bar{\epsilon}) ]$ for all $k \in [K]$, where $K$ is the total number of iterations.
\end{remark}
}

\begin{lemma} \label{lemma:Bound_y_tilde}
For all $k \in [K]$, (i) $\| \tilde{y}^k \|^2 \in
[G^{\text{\tiny L}}, G^{\text{\tiny U}}(\bar{\epsilon}) ]$, where $G^{\text{\tiny L}}$ is a small positive number, and 

\vspace{-3mm}
\begin{small}
  \begin{align}
    G^{\text{\tiny U}}(\bar{\epsilon}) := & \sum_{z \in \mathcal{Z}} \sum_{i \in C(z)} \Big\{ \big[ \max \big\{ |y^{\text{\tiny L}}_{i}|, |y^{\text{\tiny U}}_{i}| \big\} \big]^2 + \tilde{\xi}^{\text{\tiny U}}_{zi} (\bar{\epsilon})^2  + \nonumber \\
    & 2 \tilde{\xi}^{\text{\tiny U}}_{zi} (\bar{\epsilon}) \max \big\{ |y^{\text{\tiny L}}_{i}|, |y^{\text{\tiny U}}_{i}| \big\}   \Big\}, \label{UB}
  \end{align}
\end{small}  
and (ii) we have the following basic inequality:

\vspace{-3mm}
\begin{small}
\begin{align}
  & \| \lambda^{k+1}-\lambda^{\star} \|^2  \leq \| \lambda^k - \lambda^{\star} \|^2 + \alpha_k^2 \|s^k(\tilde{y}^k)\|^2 +  \nonumber \\
  &  2 \alpha_k \big( H(\lambda^k) - H(\lambda^{\star}) \big) +  2 \alpha_k \langle s^k(\tilde{y}^k) - y^k, \lambda^k-\lambda^{\star} \rangle. \label{basic_inequality} 
\end{align}
\end{small}
\end{lemma}
\proof \eqref{UB} holds from \eqref{Perturb_y}, and Remarks \ref{remark:bound_y} and \ref{remark:noise}.
\eqref{basic_inequality} holds because of the nonexpansion property of the projection and the supergradient inequality, namely, $H(\lambda) - H(\lambda^k) \leq \langle y^k , \lambda - \lambda^k \rangle$ for all $\lambda \in \Lambda$.\qed

We emphasize that $G^{\text{\tiny U}}(\bar{\epsilon})$ increases as $\bar{\epsilon}$ decreases. 

\subsubsection{Rule 1}
Under Rule 1 it follows from \eqref{basic_inequality} that

\vspace{-3mm}
\begin{small}
\begin{align}
  & \| \lambda^{k+1}-\lambda^{\star} \|^2 
  \leq \| \lambda^k - \lambda^{\star} \|^2 + \alpha_k^2 G^{\text{\tiny U}}(\bar{\epsilon}) + \nonumber \\
  & 2 \alpha_k \big( H(\lambda^k) - H(\lambda^{\star}) \big) +  2 \alpha_k \langle \tilde{\xi}^k, \lambda^k-\lambda^{\star} \rangle, \ \forall k \in \mathbb{N}.   \label{Rule1_ineq}
\end{align}
\end{small}

{\color{black}
By taking the conditional expectation on \eqref{Rule1_ineq}, one can derive the following inequality:

\vspace{-3mm}
\begin{small}
\begin{align}
  & \mathbb{E} \big[ \| \lambda^{k+1}-\lambda^{\star} \|^2 \ | \ \lambda^k \big]  \leq \| \lambda^k - \lambda^{\star} \|^2 + \alpha_k^2 G^{\text{\tiny U}}(\bar{\epsilon}) + \nonumber \\
  & 2 \alpha_k \big( {H}(\lambda^k) - {H}(\lambda^{\star}) \big), \ \forall k \in [K],
  \label{Rule1_conditional_expectation_inequality}
\end{align}
\end{small}
where the inequality holds because of $\mathbb{E}[\tilde{\xi}^k_{zi}| \lambda^k] = 0$ for all $z \in \mathcal{Z}$ and $i \in C(z)$.
By the law of total expectation with respect to $\lambda^k$ in~\eqref{Rule1_conditional_expectation_inequality}, we obtain
}

\vspace{-3mm}
\begin{small}
\begin{align}
  &\mathbb{E} \big[ \| \lambda^{k+1}-\lambda^{\star} \|^2 \big] \leq \mathbb{E} \big[ \| \lambda^k - \lambda^{\star} \|^2 \big] + \alpha_k^2 G^{\text{\tiny U}}(\bar{\epsilon}) + \nonumber \\
  & 2 \alpha_k \mathbb{E}\big[ H(\lambda^k) - H(\lambda^{\star}) \big], \ \forall k \in [K].
  \label{Rule1_expectation_inequality}
\end{align}
\end{small}
\noindent
We recursively add \eqref{Rule1_expectation_inequality} from
$k=1$ to $k=K$ to obtain

\vspace{-3mm}
\begin{small}
\begin{align}
& \mathbb{E} \big[ \| \lambda^{K+1} - \lambda^{\star} \|^2 \big] \leq \| \lambda^1 - \lambda^{\star} \|^2  +  \nonumber \\
& G^{\text{\tiny U}}(\bar{\epsilon})  \sum_{k=1}^{K} \alpha^2_k + 2 \sum_{k=1}^{K} \alpha_k \Big( \mathbb{E} \big[ H(\lambda^k) \big] - H(\lambda^{\star}) \Big). \label{Rule1_recursive_ineq}
\end{align}
\end{small}
Since $\exists \lambda^{\text{\tiny U}} : \lambda^{\text{\tiny U}} \geq \| \lambda^1 - \lambda^{\star} \|^2$ by Assumption \ref{assump:dual_bound} and $\mathbb{E} \big[ \|
\lambda^{K+1} - \lambda^{\star} \|^2 \big] \geq 0$, \eqref{Rule1_recursive_ineq} can be expressed as

\vspace{-3mm}
\begin{small}
\begin{align}
& \lambda^{\text{\tiny U}} +  G^{\text{\tiny U}}(\bar{\epsilon})  \sum_{k=1}^{K} \alpha^2_k \ \geq 2 \sum_{k=1}^{K} \alpha_k \Big({H}(\lambda^{\star}) - \mathbb{E} \big[ {H}(\lambda^k) \big] \Big) \nonumber \\
& \geq \ \Big(2 \sum_{k=1}^{K} \alpha_k \Big) \Big({H}(\lambda^{\star}) - \max_{k \in [K]} \mathbb{E}[{H}(\lambda^k)] \Big) \nonumber \\
& \geq \  \Big(2 \sum_{k=1}^{K} \alpha_k \Big) \Big({H}(\lambda^{\star}) - \mathbb{E} \big[ \max_{k \in [K]} {H}(\lambda^k) \big] \Big), \label{Rule1_chain_ineq}
\end{align}
\end{small}
\noindent
where the last inequality holds due to Jensen's inequality.
By substituting $\alpha_k = a/k$ in \eqref{Rule1_chain_ineq}, we obtain

\vspace{-3mm}
\begin{small}
\begin{align}
{H}(\lambda^{\star}) - \mathbb{E} \big[ {H}_{\text{best}}(\lambda^{K}) \big] \leq \frac{\lambda^{\text{\tiny U}} + G^{\text{\tiny U}}(\bar{\epsilon}) \sum_{k=1}^{\infty} (a/k)^2}{2 \sum_{k=1}^{K} (a/k)}, \label{Rule1_sandwich_inequality}
\end{align}
\end{small}
\noindent
where ${H}_{\text{best}}(\lambda^{K}) := \max_{k \in [K]} {H}(\lambda^k)$.
\begin{theorem} \label{thm:Rule1_convergence_1}
Algorithm \ref{algo:DPPSA} with Rule 1 provides a sequence that converges in expectation and probability, namely,

\vspace{-3mm}
\begin{small}
\begin{subequations}
\begin{align}
& \lim_{K \rightarrow \infty} \ \mathbb{E} \big[ {H}_{\text{best}} (\lambda^{K}) \big] = {H}(\lambda^{\star}), \label{Rule1_convergence_expectation}    \\
& \lim_{K \rightarrow \infty} \ \mathbb{P} \big\{  {H}(\lambda^{\star}) - {H}_{\text{best}} (\lambda^{K}) \geq \epsilon \big\} = 0, \label{Rule1_convergence_probability}
\end{align}
\end{subequations}
\end{small}
\noindent
for any $\epsilon > 0$.
Furthermore, the rate of convergence in expectation is $\mathcal{O}(G^{\text{\tiny U}}(\bar{\epsilon})/\log (K))$, where $G^{\text{\tiny U}}(\bar{\epsilon})$ increases as $\bar{\epsilon}$ decreases.
\end{theorem}
\proof See Appendix \ref{apx-thm:Rule1_convergence_1}.
\qed

To show that Algorithm \ref{algo:DPPSA} provides a sequence that converges with probability $1$, we introduce the notion of the stochastic quasi-Feyer sequence in Definition \ref{def:stochastic_quasi}.
\begin{definition}{ \text{\bf (Stochastic quasi-Feyer sequence \cite{ermoliev1988numerical})}} \label{def:stochastic_quasi}
A sequence of random vectors $\{ z^k \}_{k=1}^{\infty}$ is a stochastic quasi-Feyer sequence for a set $\mathcal{Z} \subset \mathbb{R}^n$ if $\mathbb{E}[ \|z^1\|^2] < \infty$, and for any $z \in \mathcal{Z}$, 

\vspace{-3mm}
\begin{small}
\begin{align*}
& \mathbb{E}\big[ \| z - z^{k+1} \|^2 \ | \ z^1, \ldots, z^k \big] \leq \| z - z^k \|^2 + d_k, \ \forall k \in \mathbb{N},  \\
& d_k \geq 0, \ \forall k \in \mathbb{N}, \ \ \sum_{k=1}^{\infty} \mathbb{E}[d_k] < \infty.
\end{align*}
\end{small}
\end{definition}

\begin{theorem} \label{thm:Rule1_convergence_as}
Algorithm \ref{algo:DPPSA} with Rule 1 provides a sequence that converges with probability $1$, namely,

\vspace{-3mm}
\begin{small}
\begin{align}
& \mathbb{P} \big\{ \lim_{K \rightarrow \infty} {H}_{\text{best}} (\lambda^{K}) = {H}(\lambda^{\star}) \big\} = 1. \label{Rule1_convergence_almost}
\end{align}
\end{small}
\end{theorem}
\proof See Appendix \ref{apx-thm:Rule1_convergence_as}.
\qed

\subsubsection{Rule 2} \label{sec:Rule2}
Under Rule 2 it follows from \eqref{basic_inequality} that

\vspace{-3mm}
\begin{small}
\begin{align}
  & \| \lambda^{k+1}-\lambda^{\star} \|^2  
    \leq \| \lambda^k - \lambda^{\star} \|^2 -  \frac{ \big( {H}(\lambda^{\star}) - {H}(\lambda^{k}) \big)^2}{\| \tilde{y}^k  \|^2}  + \nonumber \\
  &  2 \frac{{H}(\lambda^{\star}) - {H}(\lambda^{k})}{\| \tilde{y}^k  \|^2} \| \tilde{\xi}^{k} \| \cdot \|\lambda^k-\lambda^{\star} \| \label{Rule2_basic_ineq_1} \\
  & \leq \| \lambda^k - \lambda^{\star} \|^2 -  \frac{ \big( {H}(\lambda^{\star}) - {H}(\lambda^{k}) \big)^2}{G^{\text{\tiny U}}(\bar{\epsilon})}  +  M(\bar{\epsilon}), \nonumber
\end{align}
\end{small}
\noindent
where the first inequality holds due to the Cauchy--Schwarz inequality and the last inequality holds due to the existence of $M(\bar{\epsilon}) \in (0, \infty)$ based on Remark \ref{remark:noise}, Lemma \ref{lemma:Bound_y_tilde}, and Assumption \ref{assump:dual_bound}.
By taking the expectation and applying Jensen's inequality, we have

\vspace{-3mm}
\begin{small}
\begin{align}
& \mathbb{E} \big[ \| \lambda^{k+1} - \lambda^{\star} \|^2 \big] \leq \mathbb{E} \big[ \| \lambda^{k} - \lambda^{\star} \|^2 \big] - \nonumber \\
&    \big(  {H}(\lambda^{\star}) - \mathbb{E} \big[ {H}(\lambda^{k}) \big] \big)^2 / G^{\text{\tiny U}}(\bar{\epsilon}) + M(\bar{\epsilon}). \label{Rule2_expect_ineq}
\end{align}
\end{small}
Following the similar derivation from Rule 1, we obtain

\vspace{-3mm}
\begin{small}
\begin{align}
  {H}(\lambda^{\star}) - \mathbb{E} \big[ {H}_{\text{best}}(\lambda^{K}) \big]  \leq \sqrt{ \frac{\big(\lambda^{\text{\tiny U}}+K M(\bar{\epsilon})\big) G^{\text{\tiny U}} (\bar{\epsilon})}{K}}. \label{Rule2_sandwichi_ineq_1}
\end{align}
\end{small}
Based on \eqref{Rule2_sandwichi_ineq_1}, we state the following proposition.
\begin{proposition} \label{prop:Rule2_convergence}
Algorithm \ref{algo:DPPSA} with Rule 2 produces a sequence that converges in expectation to a point within $\sqrt{M(\bar{\epsilon})G^{\text{\tiny U}} (\bar{\epsilon})}$ of the optimal value. Since $M(\bar{\epsilon})G^{\text{\tiny U}} (\bar{\epsilon})$ increases as $\bar{\epsilon}$ decreases, it implies that there exists a trade-off between TPL and solution accuracy. 
\end{proposition}

We show, however, that the trade-off vanishes under the following assumption. 
\begin{assumption}{(Adapted from Assumption 3.1 in \cite{nedic2010effect})} \label{assump:sharp}
  There exists $\mu  > 0$ such that

  \vspace{-3mm}
  \begin{small}
  \begin{subequations}
  \begin{align}
    &\mu \|\lambda - \lambda^{\star}\| \leq {H} (\lambda^{\star}) - {H} (\lambda), \ \forall \lambda \in \Lambda, \label{linear_growth} \\
  & \| s^k(\tilde{y}^k) - {y}^k \| <  \mu/2, \ \forall k \in \mathbb{N}, \label{search_direction}
  \end{align}
  \end{subequations}
  \end{small}
  \noindent
  where the first inequality indicates that the function ${H}$ has a sharp set of maxima over a convex set $\Lambda$ and the second inequality indicates that the distance between the search direction and the supergradient is bounded.
\end{assumption}
Assumption~\ref{assump:sharp} is mild since the function $H$ is polyhedral for our case with a reasonable choice of TPL. 

\begin{theorem} \label{thm:Rule2_convergence}
Under Assumption \ref{assump:sharp}, Algorithm \ref{algo:DPPSA} with Rule 2 provides a sequence that converges in expectation, in probability, and with probability $1$.
The rate of convergence in expectation is $\mathcal{O}(G^{\text{\tiny U}}(\bar{\epsilon})/\sqrt{K})$, where $G^{\text{\tiny U}}(\bar{\epsilon})$ increases as $\bar{\epsilon}$ decreases.
\end{theorem}
\proof
See Appendix \ref{apx-thm:Rule2_convergence}.
\qed

\subsubsection{Rule 3} 
Under Rule 3 the search direction $s^k(\tilde{y}^k)$ is a linear combination of $\{\tilde{y}^{k'}\}_{k'=1}^{k-1}$.
\begin{lemma} \label{lemma:Rule3_upper_bound_s}
  Under Rule 3 we have
  \begin{align*}
  \| s^k (\tilde{y}^k) \|^2 = \|\tilde{y}^k + \zeta_k s^{k-1} (\tilde{y}^{k-1})\|^2 \leq \|\tilde{y}^k\|^2, \ \forall k \in \mathbb{N}.
  \end{align*}
  \end{lemma}
  \proof
  If $\zeta_k = 0$, then $s^k(\tilde{y}^k) = \tilde{y}^k$.
  If $\zeta_k > 0$, then

  \vspace{-3mm}
  \begin{small}
  \begin{align*}
  & \|\tilde{y}^k + \zeta_k s^{k-1} \|^2 - \|\tilde{y}^k\|^2  = \zeta_k^2 \|s^{k-1}\|^2 + 2 \zeta_k \langle  s^{k-1}, \tilde{y}^k \rangle \\
  &= \chi_k^2 \frac{\langle s^{k-1}, \tilde{y}^k \rangle ^2}{\|s^{k-1}\|^2} - 2 \chi_k \frac{\langle s^{k-1}, \tilde{y}^k \rangle ^2}{\|s^{k-1}\|^2} \leq 0 ,
  \end{align*}
  \end{small}
  \noindent
  where the last inequality holds since $\chi_k \in [0,2]$ as defined in Table \ref{tab:three_rules}.
  \qed

From Lemma \ref{lemma:Rule3_upper_bound_s}, similar results from Rule 2 can be derived.
Under Rule 3 it follows from \eqref{basic_inequality} that

\vspace{-3mm}
\begin{small}
\begin{align}
  & \| \lambda^{k+1}-\lambda^{\star} \|^2  \leq  \| \lambda^k - \lambda^{\star} \|^2 -  \frac{ \big( {H}(\lambda^{\star}) - {H}(\lambda^{k}) \big)^2}{\| s^k (\tilde{y}^k)  \|^2}  + \nonumber \\
  & 2 \frac{{H}(\lambda^{\star}) - {H}(\lambda^{k})}{\| s^k(\tilde{y}^k)  \|^2} \| s^k(\tilde{y}^k)- {y}^k \| \cdot \|\lambda^k-\lambda^{\star} \| \label{Rule3_basic_ineq_1} \\
  & \leq \| \lambda^k - \lambda^{\star} \|^2 -  \frac{ \big( {H}(\lambda^{\star}) - {H}(\lambda^{k}) \big)^2}{G^{\text{\tiny U}}(\bar{\epsilon})}  +  R(\bar{\epsilon}), \nonumber
\end{align}
\end{small}
\noindent
where the last inequality holds due to  Lemma \ref{lemma:Rule3_upper_bound_s} and  the existence of $R(\bar{\epsilon}) \in (0, \infty)$ based on the boundness of $s^k(\tilde{y}^k)$ by its construction, Lemma \ref{lemma:Bound_y_tilde}, and Assumption \ref{assump:dual_bound}.
We emphasize that  \eqref{Rule3_basic_ineq_1} is similar to \eqref{Rule2_basic_ineq_1}. Thus one can derive results similar to \eqref{Rule2_sandwichi_ineq_1}, Proposition \ref{prop:Rule2_convergence}, and Theorem \ref{thm:Rule2_convergence} under Rule 3.

\begin{remark} \label{remark:ACOPF_convergence}
We remark that all the results related to the convergence of DP-PS also hold when solving AC OPF described in Remark \ref{remark:ACOPF_model}. However, the strong duality does not hold for AC OPF, so the consensus constraints may not be satisfied at termination.
\end{remark}

\vspace{-5mm}
\section{Numerical Experiments} \label{sec:numerical}
To support our findings from Sections \ref{sec:distributed_algo} and \ref{sec:algorithm}, 
{\color{black} we showcase that increasing TPL of DP-PS does not affect the solution accuracy, although it does affect computation.}
In all the experiments, we solve optimization models by IPOPT \cite{wachter2006implementation} via Julia 1.5.0 on a personal laptop with an Intel Core i9 CPU and 64 GB of RAM.

\subsubsection{Experimental Settings} \label{sec:experimental_setting}

For the power network instances, we consider case 14 and case 118 from Matpower \cite{zimmerman2010matpower}.
The optimal objective values of SOC OPF~\eqref{ACOPF-rect} are obtained by utilizing IPOPT: $Z^* = 8075.1$ for case 14 and $Z^* = 129341.9$ for case 118.
The networks are decomposed as described in Table \ref{tab:decomposition}.

\begin{table}[htpb]
  \centering
  \caption{Set $\mathcal{N}_z$ of buses for each zone $z \in \mathcal{Z}$.}
  \label{tab:decomposition}
  \begin{tabular}{|l|l|l|}
  \hline
         & case 14 & case 118 \cite{guo2015intelligent} \\    \hline
  Zone 1 & \{1--5\} & \{1--33\}, \{113--115\},\{117\} \\
  Zone 2 & \{7--10\} & \{34--75\}, \{116\}, \{118\}\\
  Zone 3 & \{6\}, \{11--14\} & \{76--112\}  \\
  \hline
  \end{tabular}
\end{table}

We consider an adversary who aims to estimate $\bar{D}_{14} = 47.8$ MW for case 14 and $\bar{D}_{13} = 39$ MW for case 118, respectively, by solving the adversarial problem \eqref{Adversarial} for $|\widehat{\mathcal{K}}(T)|$ times, where $T$ is any integer number less than total iterations $K$ of DP-PS and 

\vspace{-3mm}
\begin{small}
\begin{align}
\widehat{\mathcal{K}}(T) :=  \cup_{t=1}^{\lfloor K/T \rfloor} \{ (t-1)T+1, \ldots, tT  \}. \label{def_K}
\end{align}
\end{small}
Recall that $\widehat{\mathcal{K}}$ in Section \ref{sec:motivating} is $\widehat{\mathcal{K}}(1)$.

{\color{black}
On the other hand, we aim to protect demand data from the adversary by using the proposed DP-PS.
First, we compute $\bar{\Delta}^k_{zi}$ in \eqref{model:Sensitivity_problem} with $\beta=5\%$. 
Second, we consider various $\bar{\epsilon} \in \{0.01, 0.05, 0.1, 1, 10, \infty\}$ of DP-PS, where $\bar{\epsilon} = \infty$ represents a non-private PS and smaller $\bar{\epsilon}$ ensures stronger data privacy.
Note that DP-PS with $K$ total iterations and $1/\bar{\epsilon}$ TPL provides $\bar{\epsilon}$-DP (resp., $K\bar{\epsilon}$-DP) against an adversary with $T=1$ (resp., $T=K$) in \eqref{def_K}.
We use Rule 3 depicted in Section \ref{sec:algorithm} for our experiments.
}


{\color{black}
\subsubsection{Comparison with DP-ADMM}
We compare the proposed DP-PS with the existing DP-ADMM \cite{dvorkin2019differentially} (see Appendix \ref{apx-DP-ADMM} for details on DP-ADMM).

In Figure \ref{fig:DPADMM_Comparison} we report the objective values resulting by DP-PS and DP-ADMM for solving the case-14 and case-118 instances.
When $\bar{\epsilon} = 0.01$ (i.e., larger noises are introduced for stronger data privacy), the objective value of DP-ADMM significantly fluctuates and is even larger than the optimal objective value $Z^*$ of the SOC OPF model.
This implies that the sequence provided by DP-ADMM does not converge especially when stronger data privacy is required.
In contrast, the proposed DP-PS always provides a lower bound on $Z^*$ and the sequence converges to $Z^*$.

\begin{figure}[!h]  
	\centering
	\begin{subfigure}[b]{0.24\textwidth}
		\centering
		\includegraphics[width=\textwidth]{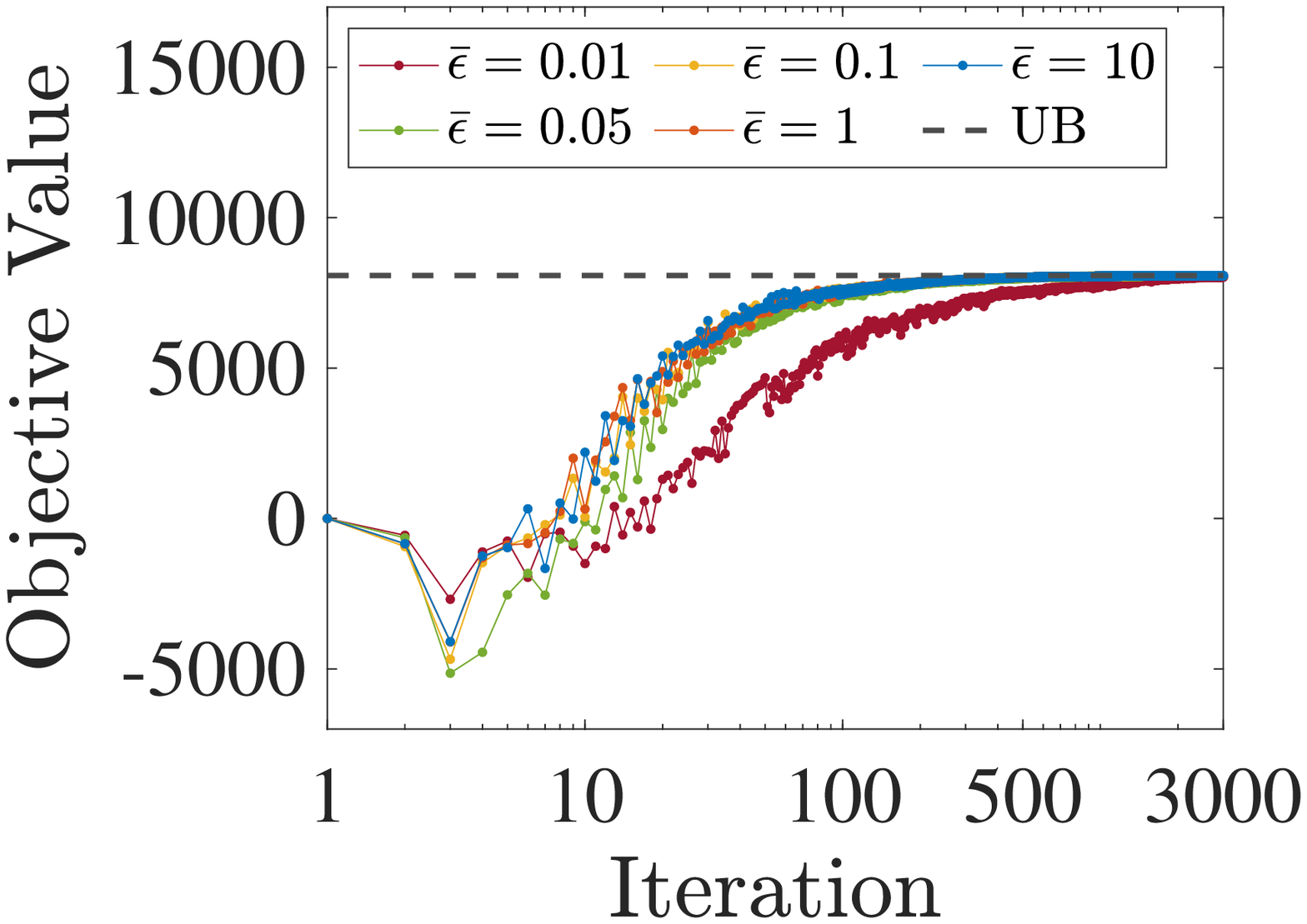}
		\includegraphics[width=\textwidth]{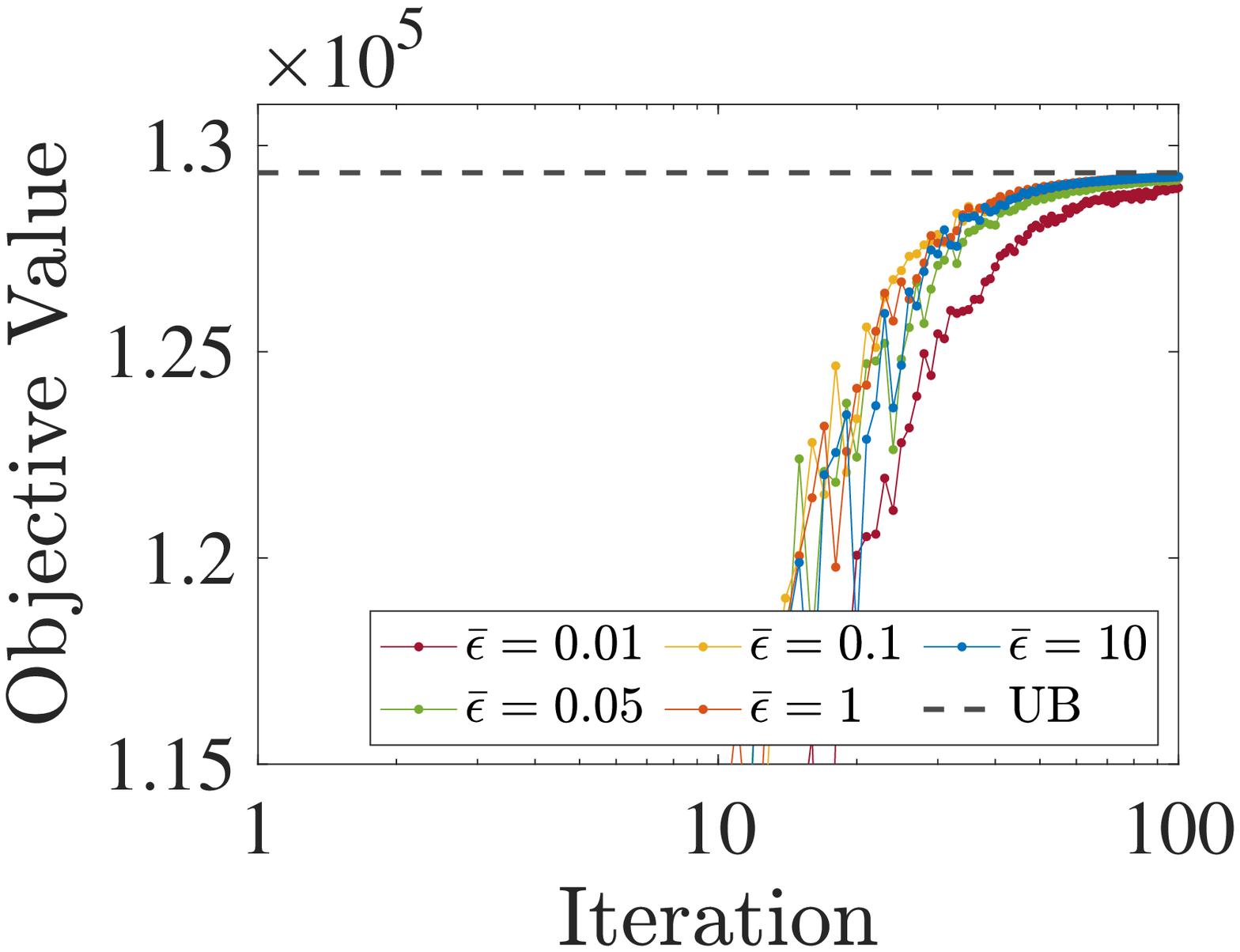}		
		\caption{DP-PS}
	\end{subfigure}
	\begin{subfigure}[b]{0.24\textwidth}
		\centering
		\includegraphics[width=\textwidth]{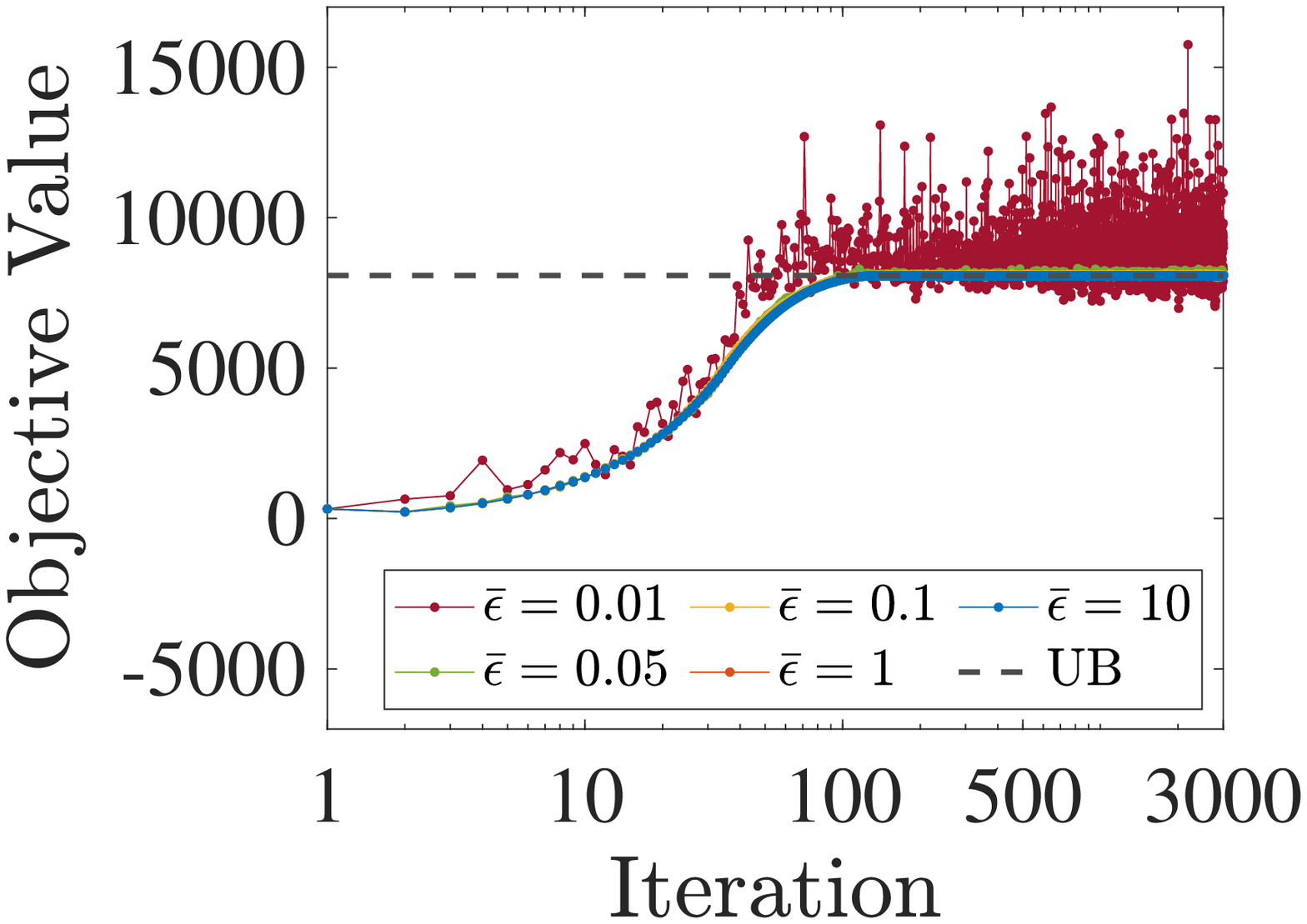}		
		\includegraphics[width=\textwidth]{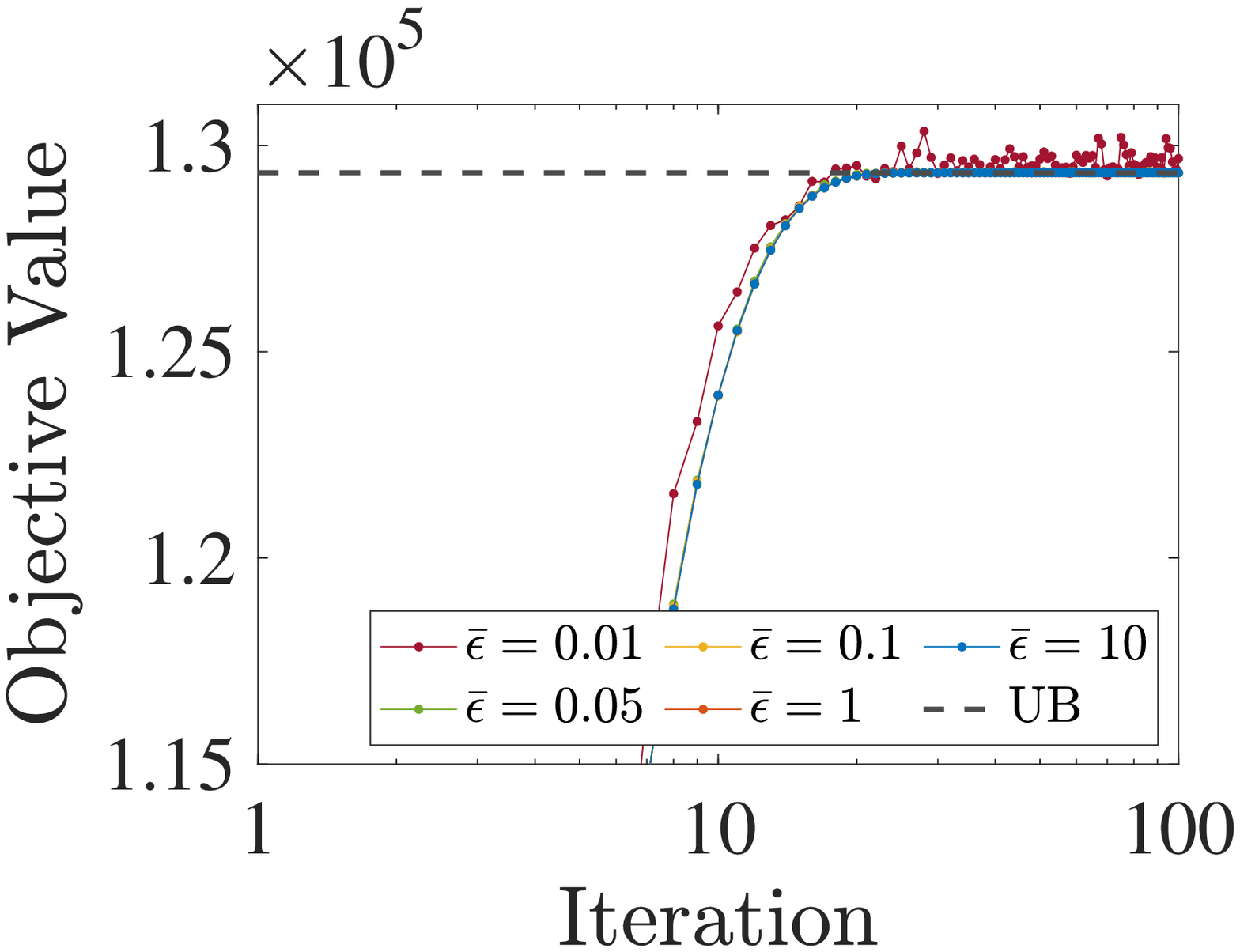}
		\caption{DP-ADMM}
	\end{subfigure}  
	\caption{Objective values computed by DP-PS (left) and DP-ADMM (right) that solve case 14 (top) and case 118 (bottom) under various $\bar{\epsilon}$.}	
  \label{fig:DPADMM_Comparison}
  \end{figure}
}
\subsubsection{Convergence of DP-PS}
We demonstrate the numerical support for Theorem \ref{thm:Rule2_convergence} that increasing TPL does not affect the solution accuracy of DP-PS, although it does affect computation.

\begin{figure}[htpb]
  \centering  
  \includegraphics[scale=0.23]{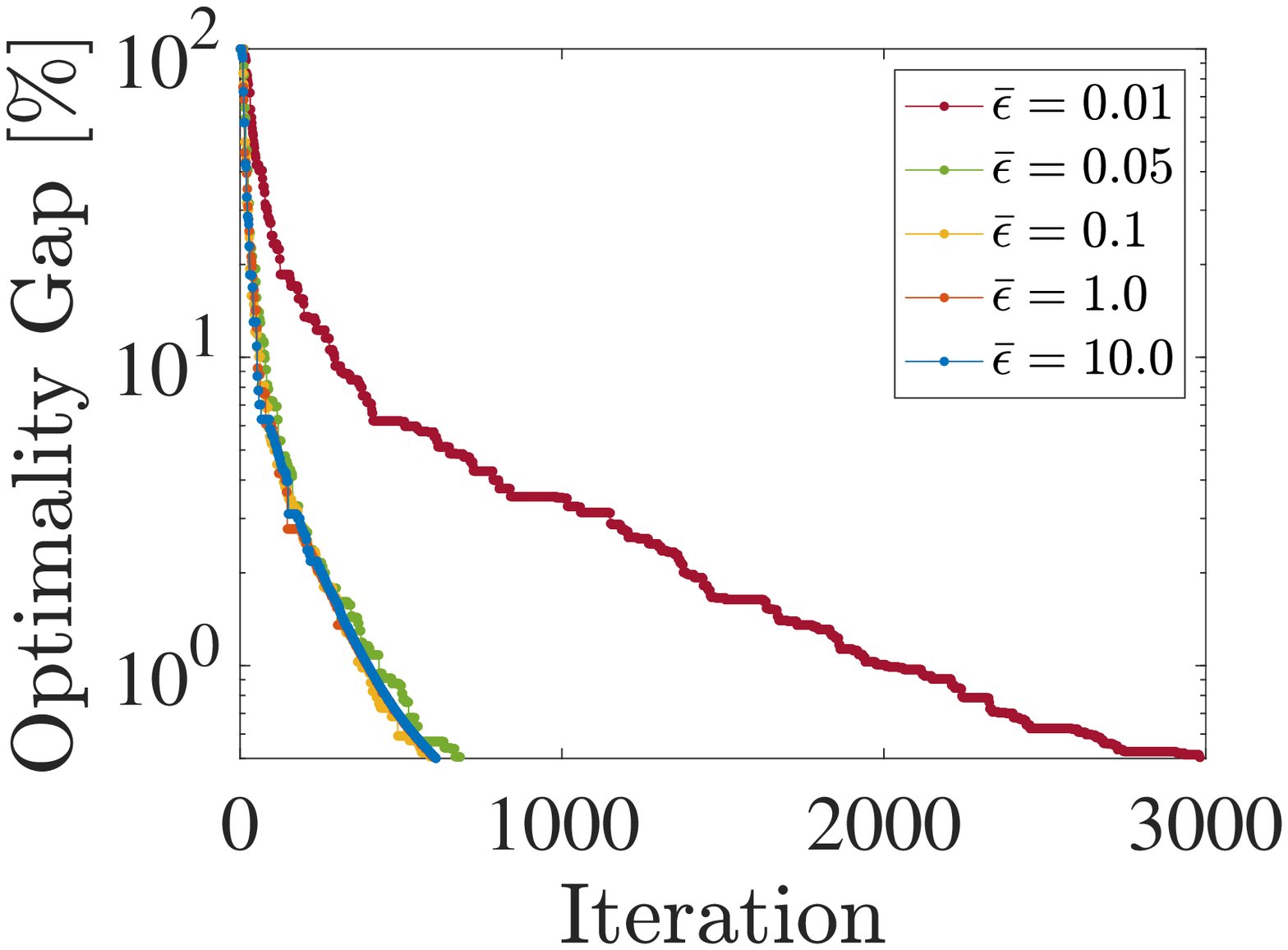}
	\includegraphics[scale=0.23]{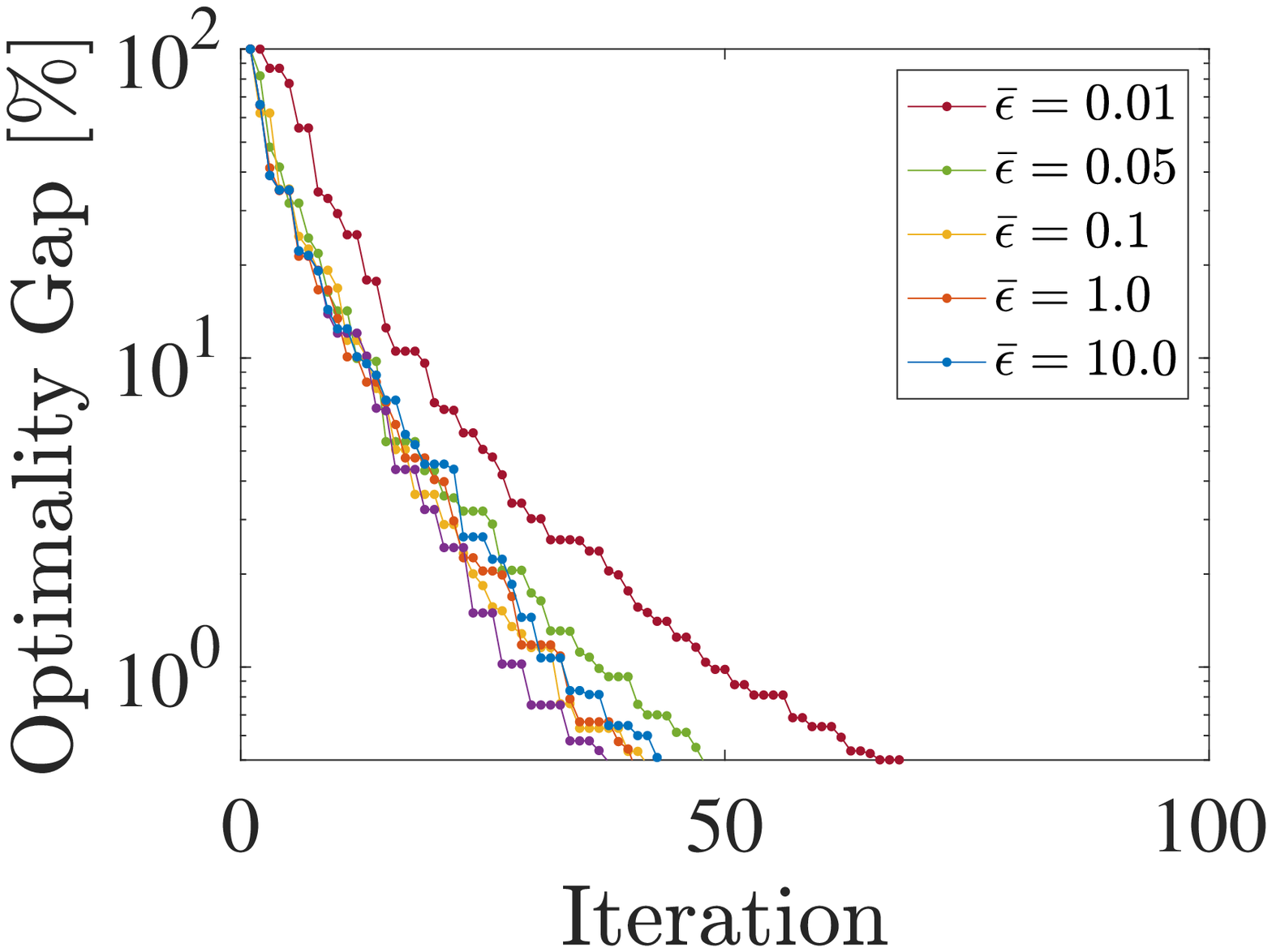}
  \caption{Optimality gap of DP-PS that solves case 14 (left) and case 118 (right) under various $\bar{\epsilon}$.}
  \label{fig:approximation_error}
\end{figure}

\noindent \textbf{Solution Accuracy:}
In Figure \ref{fig:approximation_error} we report the optimality gap at each iteration $k$ of DP-PS.
The results show that the sequence generated by DP-PS converges regardless of the {\color{black} $\bar{\epsilon}$ value}.
{\color{black}
We also discuss the impact of the number of zones on the convergence of DP-PS in Section \ref{apx-zones}.
}

\begin{figure}[htpb]
\centering
\includegraphics[scale=0.23]{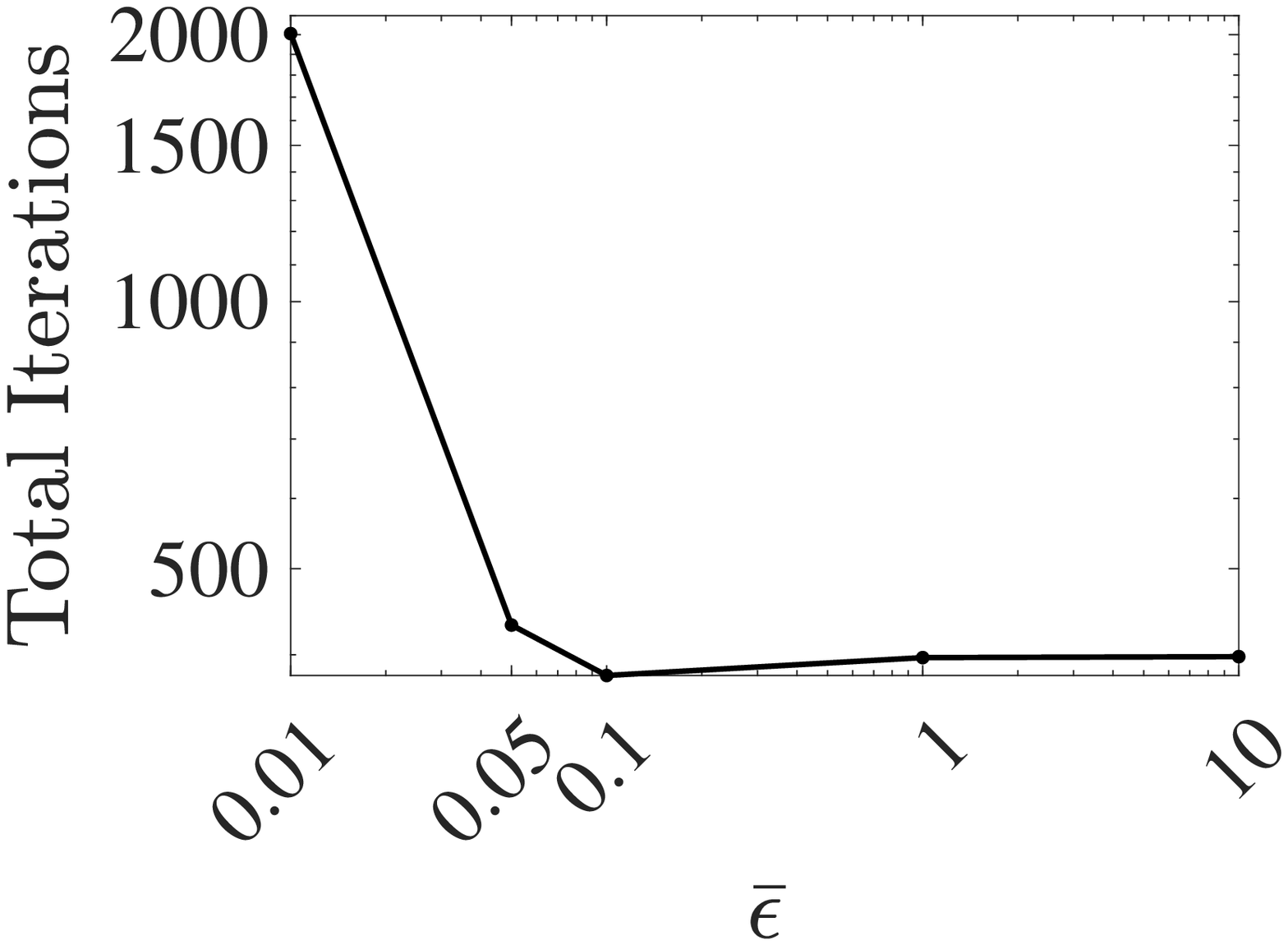}
\includegraphics[scale=0.23]{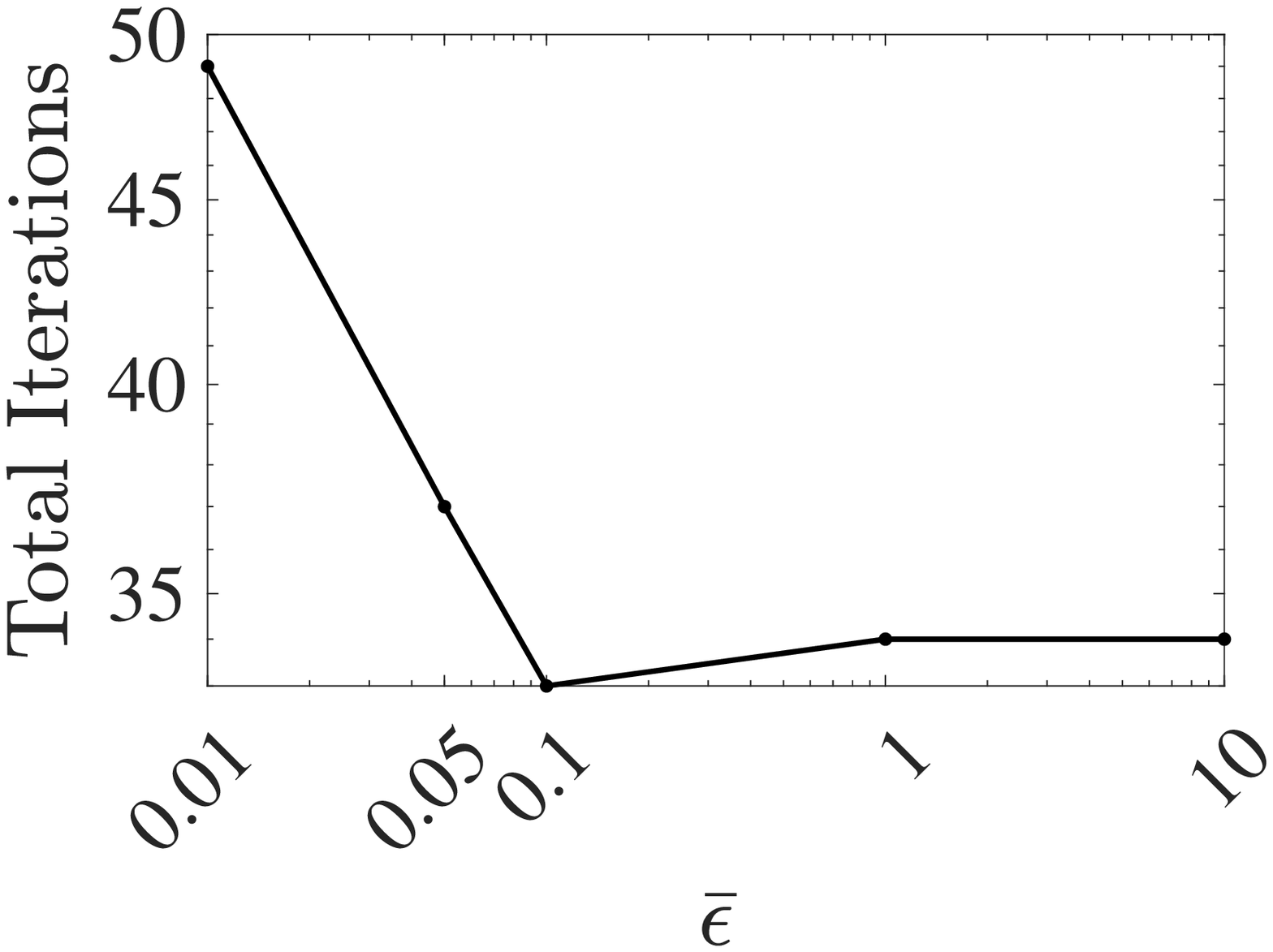}
\caption{Total iterations of DP-PS for solving case 14 (left) and case 118 (right) required to enter within $1\%$ of the optimality gap.  }
\label{fig:Iterations}
\end{figure}

\noindent \textbf{Computation:} 
Figure~\ref{fig:approximation_error} demonstrates that DP-PS with {\color{black}smaller $\bar{\epsilon}$} requires more iterations to converge. 
In Figure \ref{fig:Iterations} we report the total number of iterations required for DP-PS to converge to a solution within $1\%$ of the optimality gap. 
{\color{black}The results show the decreasing trends of total iterations as $\bar{\epsilon}$ increases.}
This implies that there exists a trade-off between TPL and computation.

\subsubsection{Data Privacy Preservation} 
We numerically show that increasing TPL provides higher data privacy.

\begin{figure}[htpb]
  \centering  
  \includegraphics[scale=0.23]{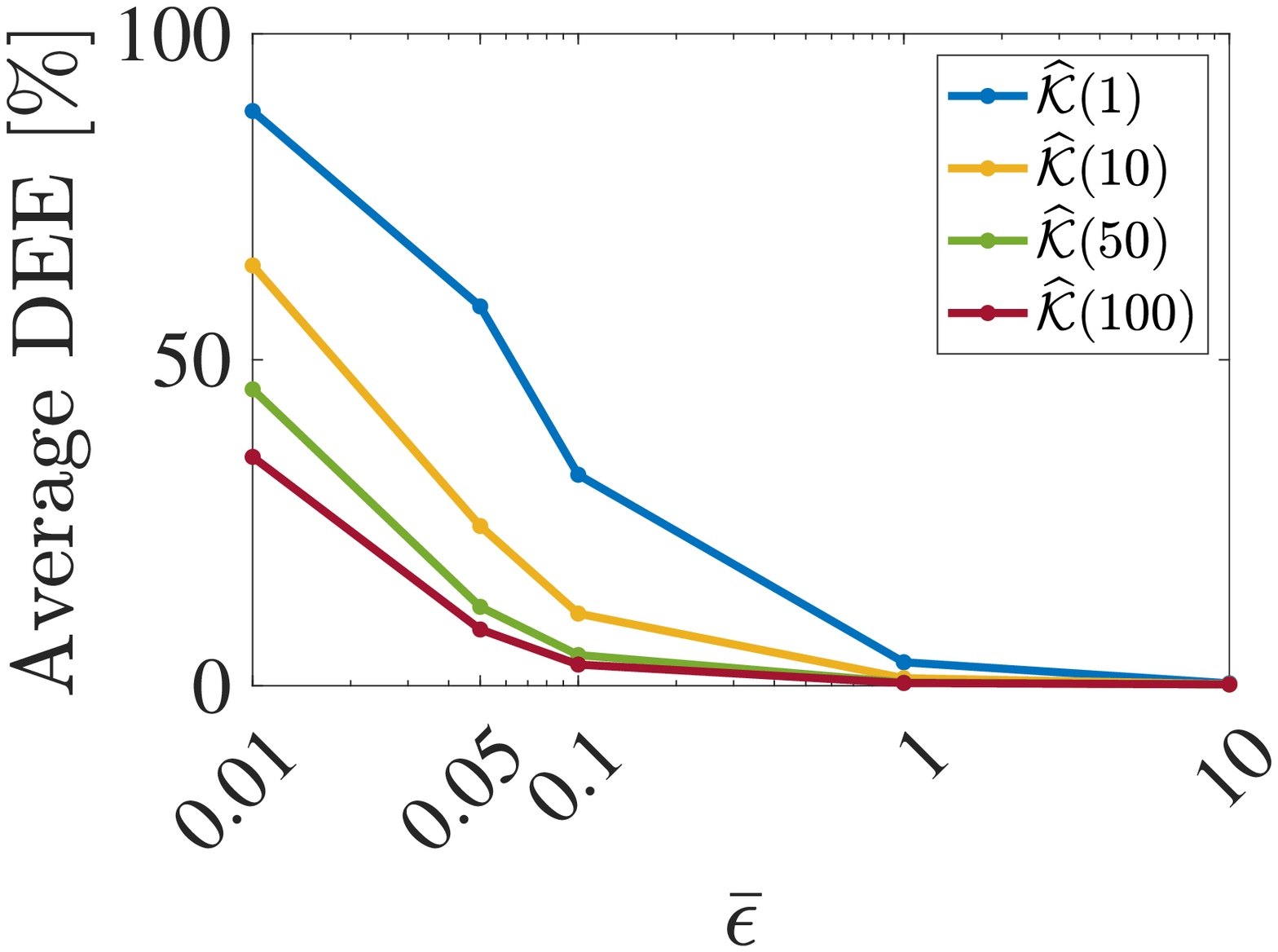}
  \includegraphics[scale=0.23]{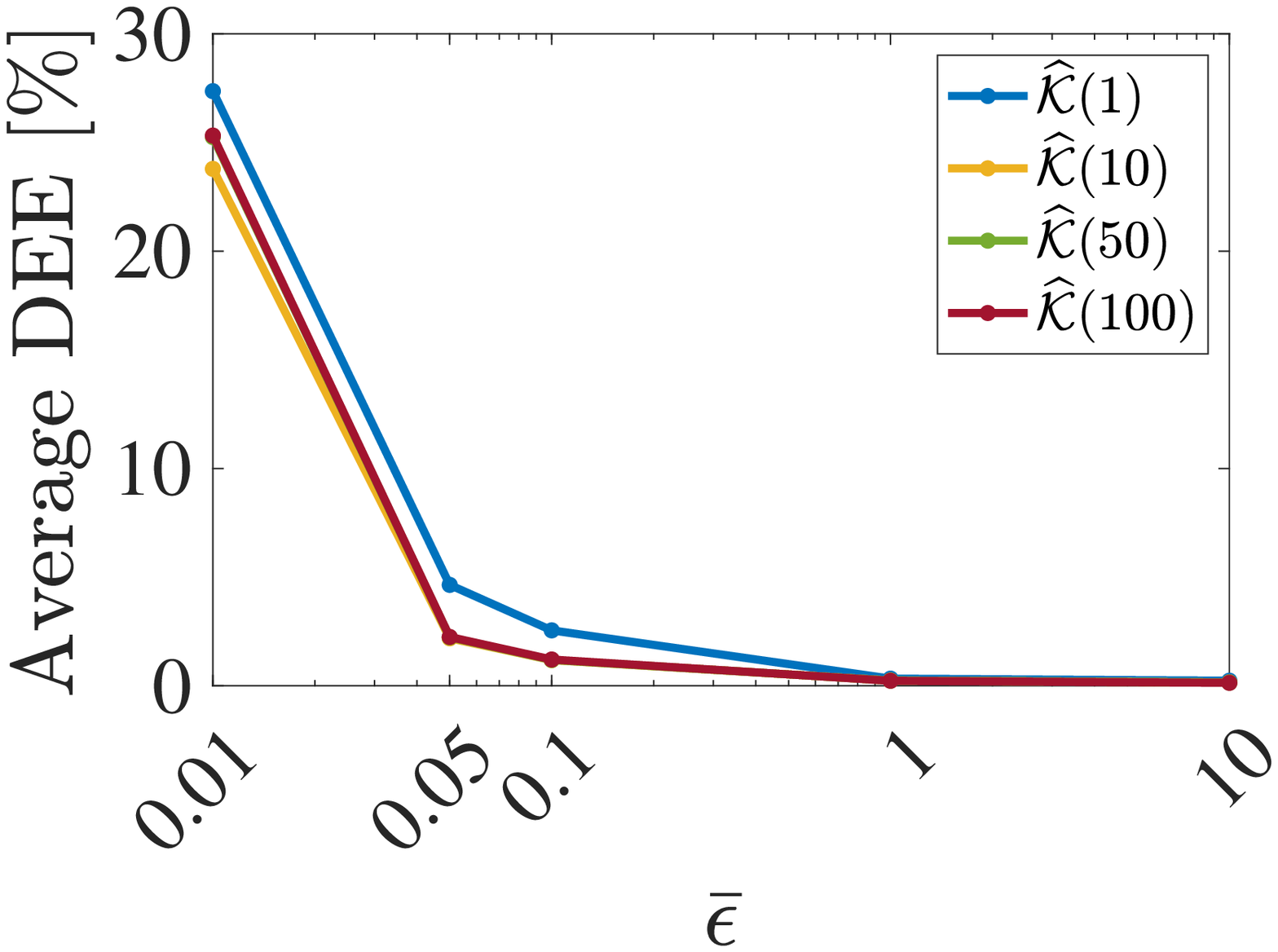}
  \caption{Average demand estimation error under various $\widehat{\mathcal{K}}(T)$ for case 14 (left) and case 118 (right).}
  \label{fig:various_K}
\end{figure}

First, we consider various $\widehat{\mathcal{K}}(T)$ when constructing the adversarial problem \eqref{Adversarial}.
As $T$ increases, theoretically, the accuracy of the demand estimated by solving \eqref{Adversarial} with $\mathcal{K} \in \widehat{\mathcal{K}}(T)$ increases.
We report in Figure \ref{fig:various_K} an average demand estimation error (DEE): $\sum_{\mathcal{K} \in \widehat{\mathcal{K}}(T)} \text{DE}(\mathcal{K}) / |\widehat{\mathcal{K}}(T)|$, where $\text{DE}(\mathcal{K})$ is defined in \eqref{DE}.
The results show  
(i) decreasing trends of the average DEE as $\bar{\epsilon}$ increases for fixed $T$ and
(ii) decreasing trends of the average DEE as $T$ increases for fixed $\bar{\epsilon}$.
Moreover, the average DEE for fixed $\bar{\epsilon}$ seems to converge to a point as $T$ increases.
{\color{black}The results imply that increasing TPL produces stronger data privacy (e.g., see $\widehat{\mathcal{K}}(100)$ in Figure \ref{fig:various_K} (right) when $\bar{\epsilon}=0.01$).}

 
\subsubsection{Summary}
\begin{figure}[htpb]
  \centering  
  \includegraphics[scale=0.23]{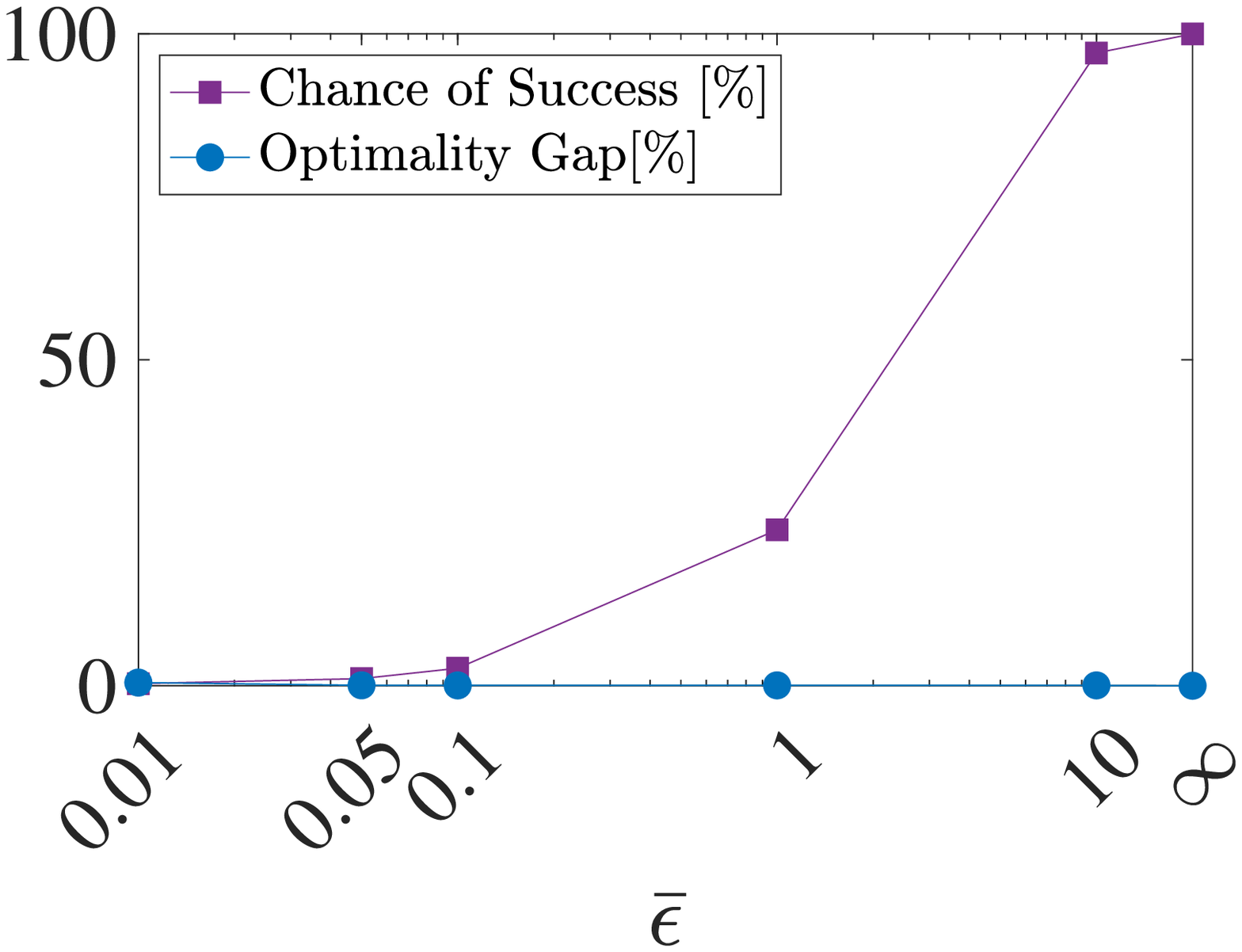}
  \includegraphics[scale=0.23]{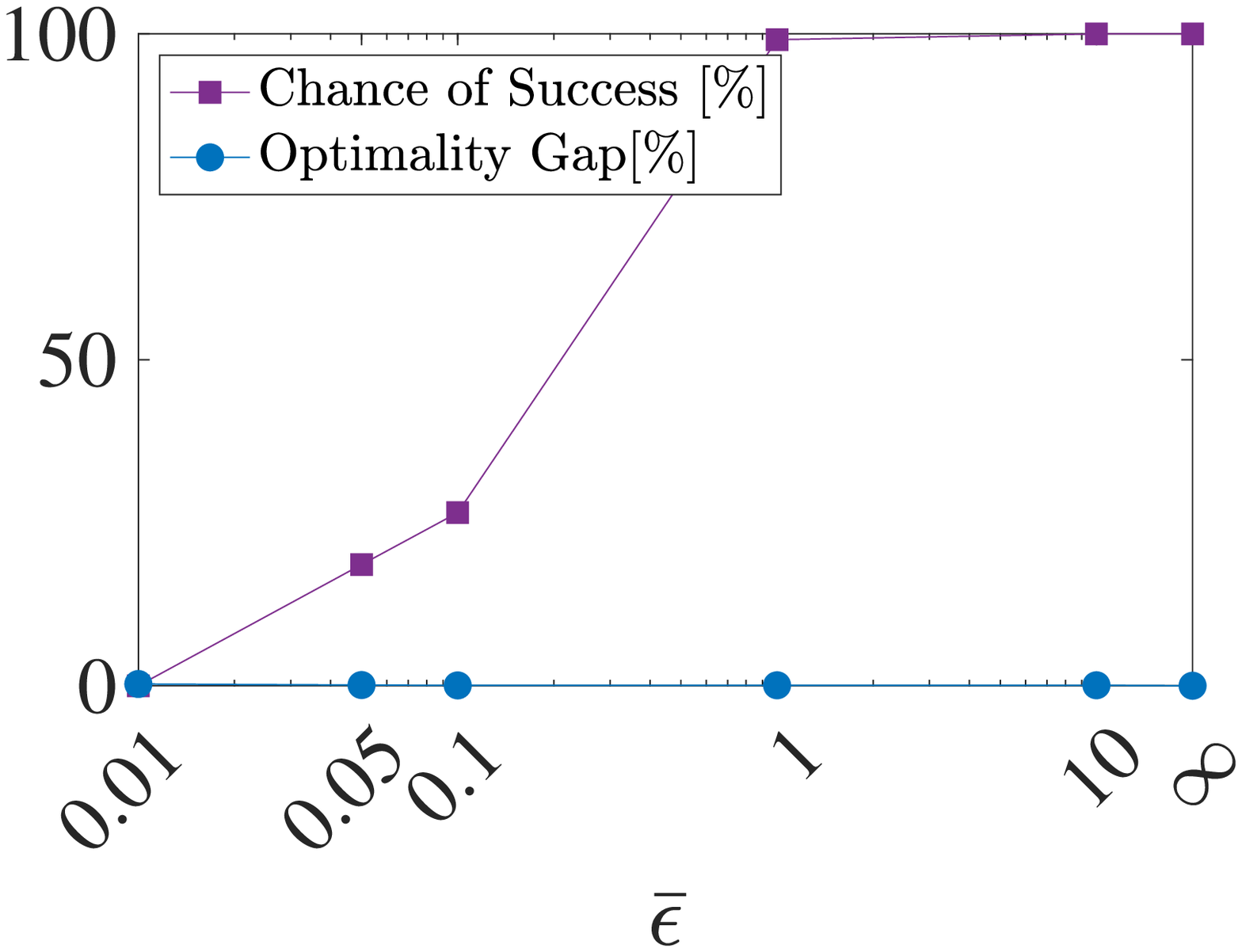}
  \caption{Summary for case 14 (left) and case 118 (right).}
  \label{fig:summary}
\end{figure}

We report in Figure \ref{fig:summary} the optimality gap at the termination of DP-PS and the adversarial's chance of success (CoS) defined as follow:

\vspace{-3mm}
\begin{small}
\begin{align*}
\text{CoS}(\overline{G}) = 100 \times \sum_{T \in \mathcal{T}} \sum_{\mathcal{K} \in \widehat{\mathcal{K}}(T)} \mathcal{I}(\text{DE}(\mathcal{K})) \leq \overline{G}) / \sum_{T \in \mathcal{T}} |\widehat{\mathcal{K}}(T)|,
\end{align*}
\end{small}
\noindent
where
$\overline{G}$ is a prespecified value (e.g., $\overline{G}=1\%$),
$\mathcal{T}$ is a collection of various $T$, $\widehat{\mathcal{K}}(T)$ is in \eqref{def_K}, $\text{DE}(\mathcal{K})$ is in \eqref{DE}, $\mathcal{I}(\text{DE}(\mathcal{K}) \leq \overline{G}) = 1$ if $\text{DE}(\mathcal{K}) \leq \overline{G}$ and $\mathcal{I}(\text{DE}(\mathcal{K}) \leq \overline{G}) = 0$ otherwise.
The results demonstrate that as $\bar{\epsilon}$ decreases, the adversarial's chance of successful demand estimation decreases while the optimality gap still remains the same. 

{\color{black}
\subsubsection{The impact of the number of zones on the convergence} \label{apx-zones}
As the number of zones increases, the convergence of DP-PS can become slower. 
To see this, we use a partitioning algorithm \cite{karypis1998multilevelk} to generate different zones of the power systems. This algorithm is available at \texttt{Metis.jl}.
In Figure \ref{fig:Zones}, we report the optimality gap produced by DP-PS with $\bar{\epsilon}=0.1$ for case 14 and case 118 instances decomposed by $|\mathcal{Z}| \in \{3,5,7\}$ and $|\mathcal{Z}| \in \{3,10,50\}$, respectively.
Note that there are about 2 buses for each zone of case 14 and case 118 systems when $|\mathcal{Z}|=7$ and $|\mathcal{Z}|=50$, respectively.
We observe that the optimality gap slightly increases as the number of zones increases, but not significantly. 

\begin{figure}[htpb]
  \centering
  \includegraphics[scale=0.23]{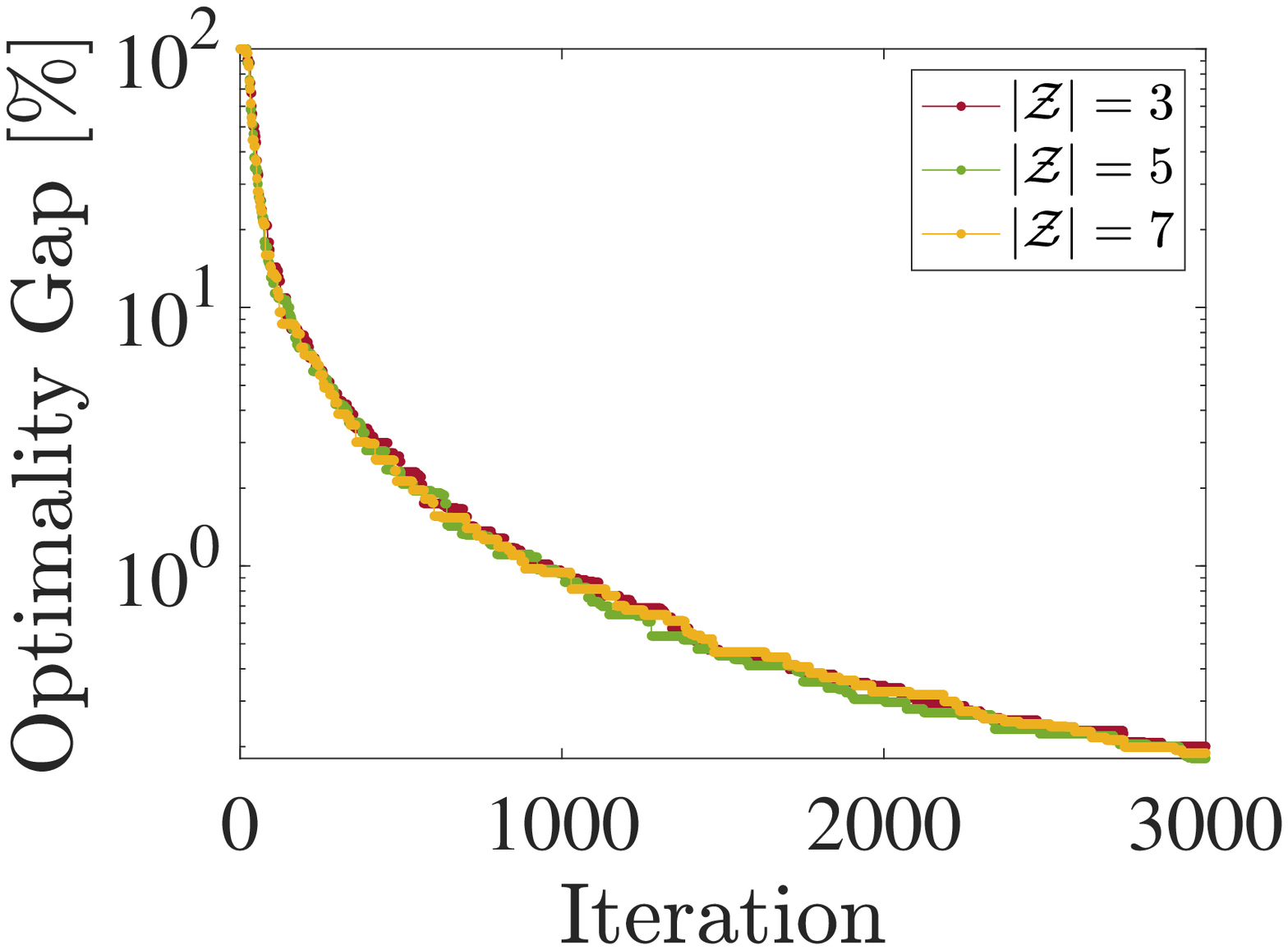}
	\includegraphics[scale=0.23]{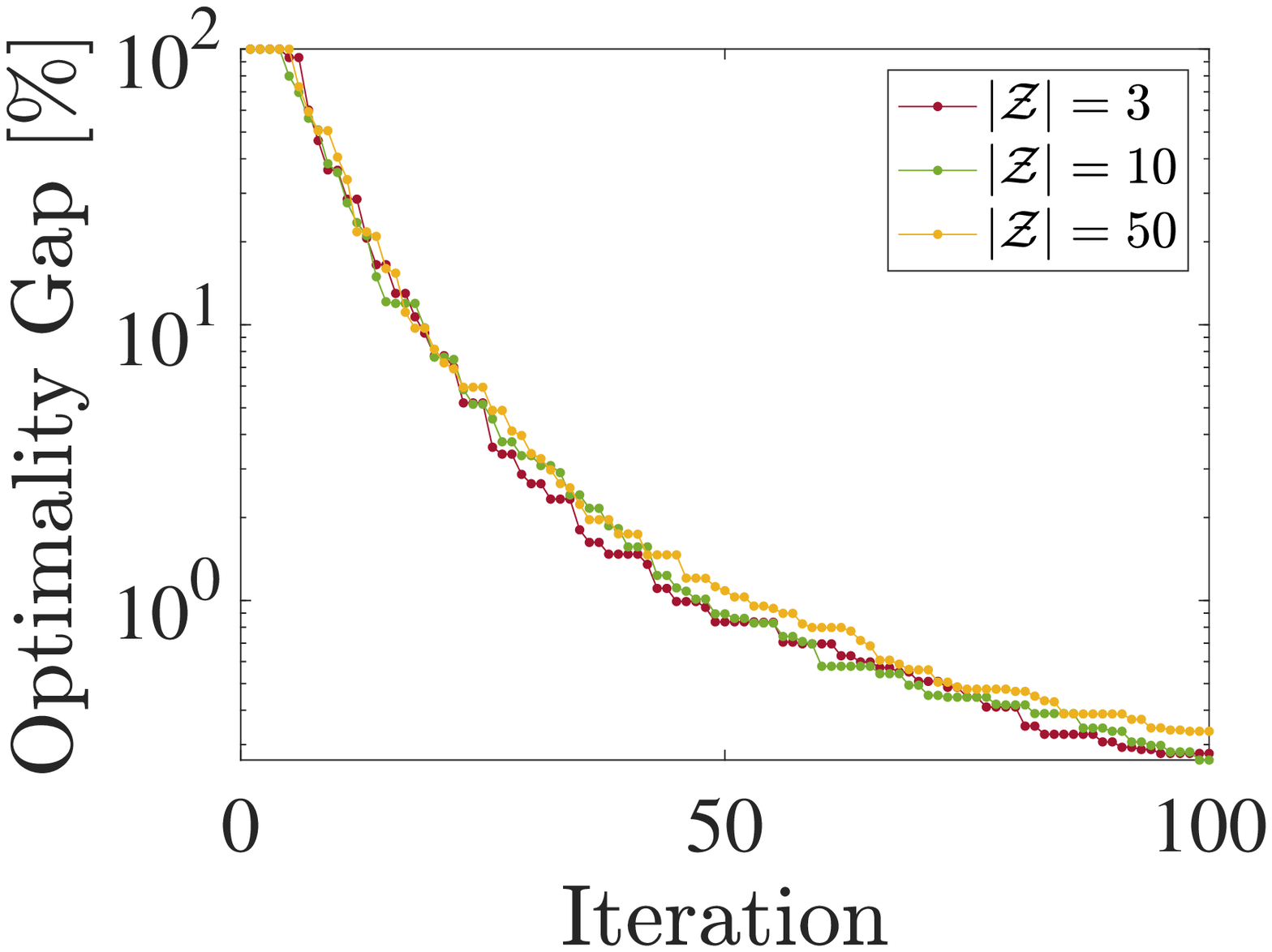}
  \caption{The impact of the number of zones on the convergence when $\bar{\epsilon}=0.1$ for case 14 (left) and case 118 (right).}
  \label{fig:Zones}
\end{figure}

}

\subsubsection{AC OPF}
\begin{figure}[htpb]
  \centering
  \includegraphics[scale=0.23]{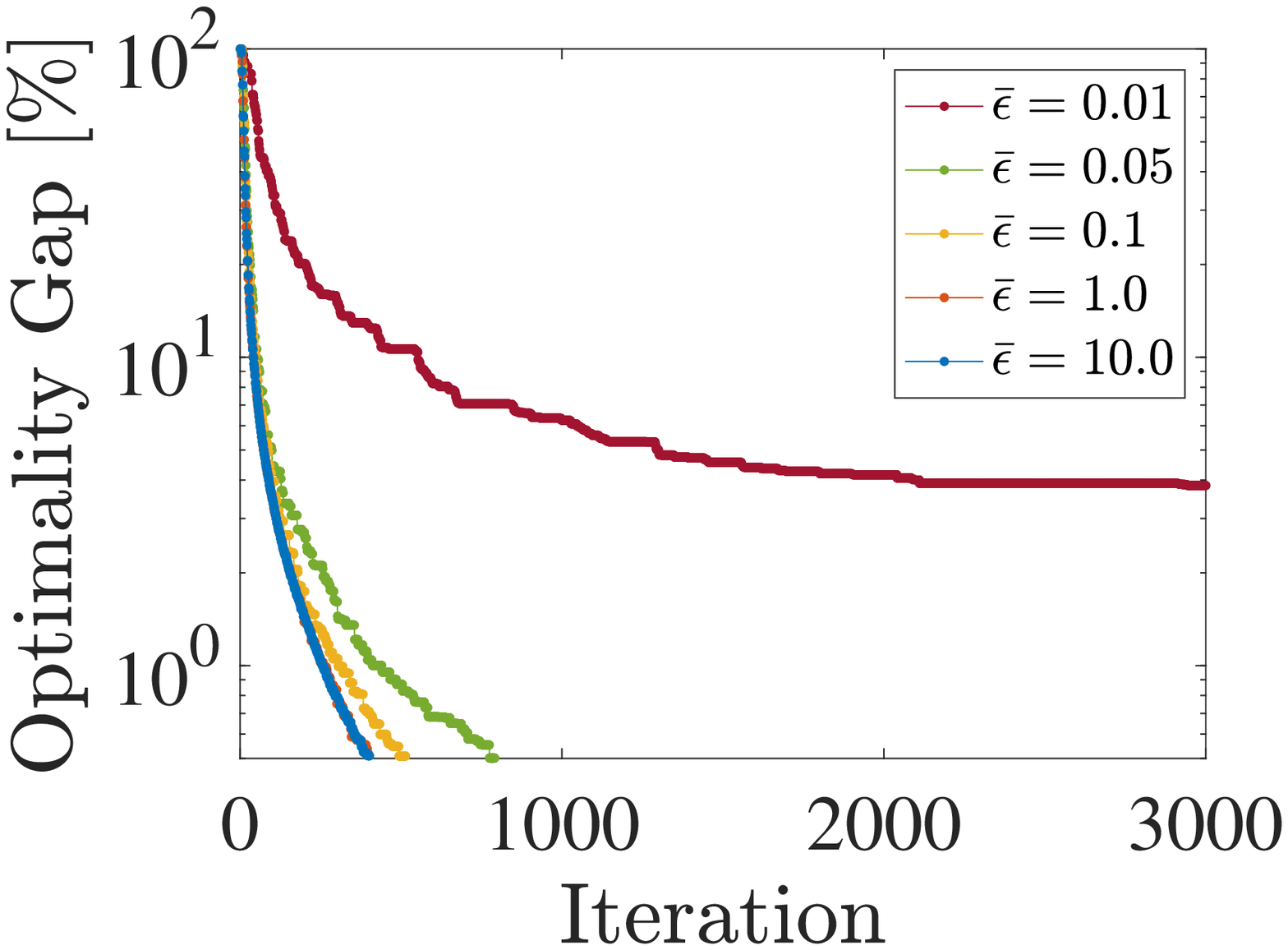}
	\includegraphics[scale=0.23]{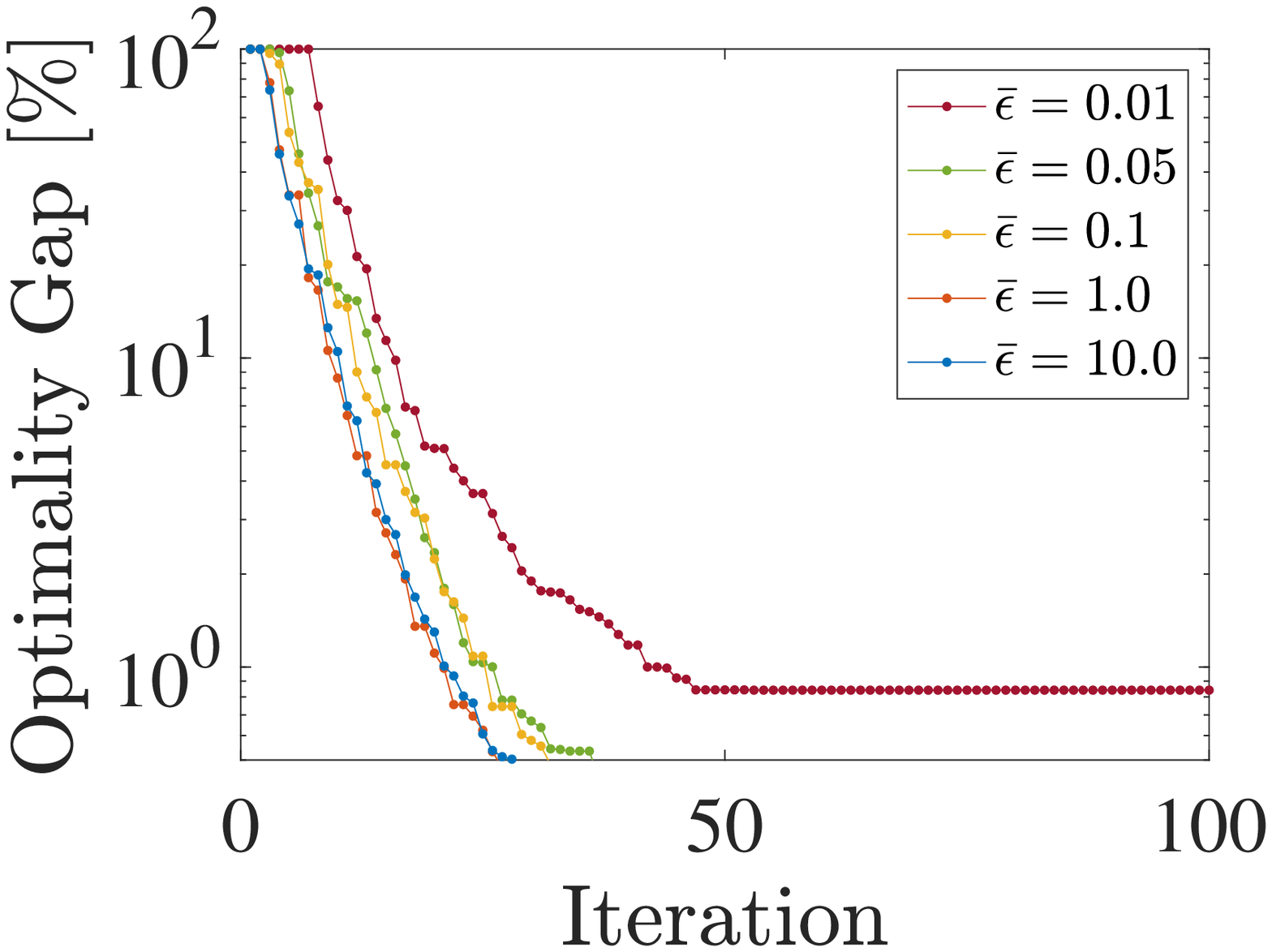}
  \caption{Approximation error of DP-PS that solves case 14 (left) and case 118 (right) under various $\bar{\epsilon}$.}
  \label{fig:AC_Approximation_Error}
\end{figure}

\begin{figure}[htpb]
  \centering  
  \includegraphics[scale=0.23]{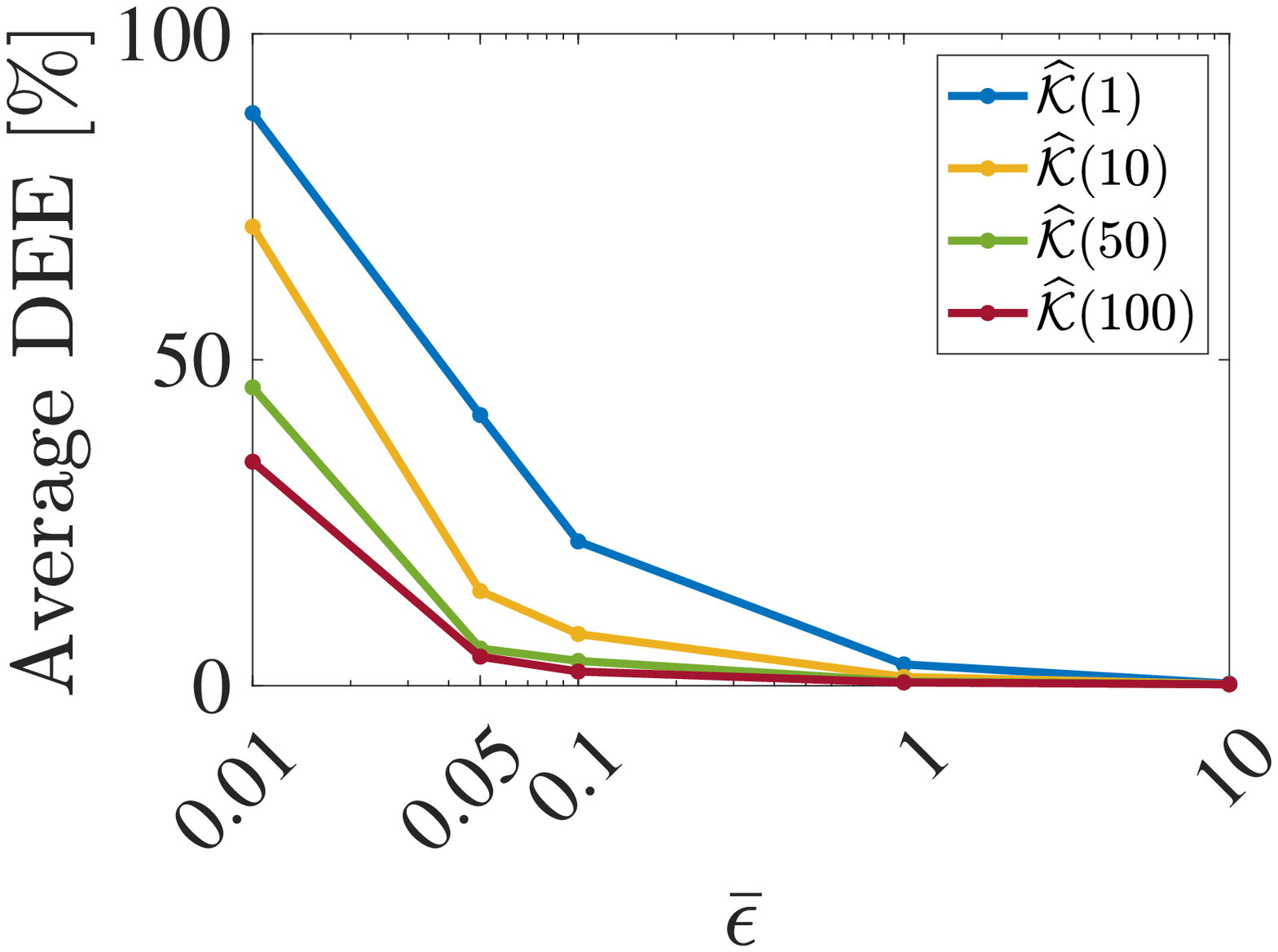}
	\includegraphics[scale=0.23]{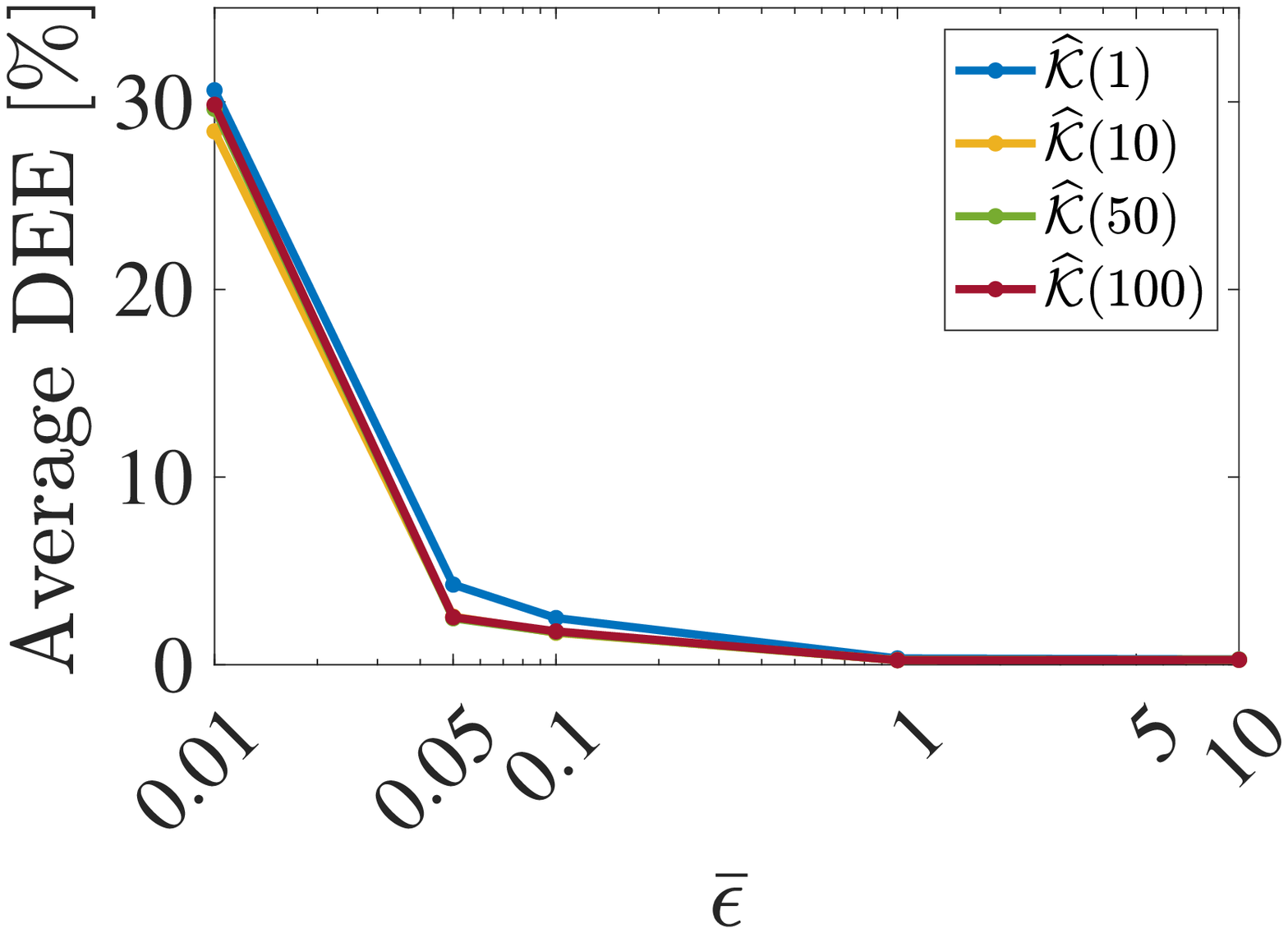}
  \caption{Average demand estimation error under various $\widehat{\mathcal{K}}(T)$ for case 14 (left) and case 118 (right).}
  \label{fig:AC_various_K}
\end{figure}

In this section we show that the convergence and the data privacy preservation of DP-PS are also achieved when solving AC OPF (see Remarks \ref{remark:ACOPF_model} and \ref{remark:ACOPF_convergence}).
To this end we present Figures \ref{fig:AC_Approximation_Error} and \ref{fig:AC_various_K}, which are counterparts of Figures \ref{fig:approximation_error} and \ref{fig:various_K}, respectively.
{\color{black} When $\bar{\epsilon}=0.01$, more iterations may be required for the convergence. Again, the sequence produced by DP-PS may converge to an infeasible point for AC OPF.}


\section{Conclusion} \label{sec:conclusion}
We studied a privacy-preserving distributed OPF and proposed a differentially private projected subgradient (DP-PS) algorithm that includes a solution encryption step.
In this algorithm Laplacian noise is introduced to encrypt solution exchanged within the algorithm, which leads to $\bar{\epsilon}$-differential privacy on data.
The target privacy level of DP-PS is chosen by users, which affects not only the data privacy but also the convergence of the algorithm.
{\color{black}
We showed that a sequence provided by DP-PS converges to an optimal solution regardless of the $\bar{\epsilon}$ value, but more iterations are required for the convergence as $\bar{\epsilon}$ decreases.
Also, we demonstrated that, as $\bar{\epsilon}$ decreases, the adversarial's chance of successful demand data inference decreases while the optimality gap remains the same.
Finally, as stated in Remarks \ref{remark:ACOPF_model} and \ref{remark:ACOPF_convergence}, the proposed DP-PS can lead to an infeasible solution with respect to the AC OPF model, calling for privacy-preserving algorithms for the model.
}


\bibliographystyle{IEEEtran}
\bibliography{References.bib}

\vspace{0.5in}
\noindent\fbox{\parbox{0.47\textwidth}{
The submitted manuscript has been created by UChicago Argonne, LLC, Operator of Argonne National Laboratory (``Argonne''). Argonne, a U.S. Department of Energy Office of Science laboratory, is operated under Contract No. DE-AC02-06CH11357. The U.S. Government retains for itself, and others acting on its behalf, a paid-up nonexclusive, irrevocable worldwide license in said article to reproduce, prepare derivative works, distribute copies to the public, and perform publicly and display publicly, by or on behalf of the Government. The Department of Energy will provide public access to these results of federally sponsored research in accordance with the DOE Public Access Plan (http://energy.gov/downloads/doe-public-access-plan).}
}

\appendices

\section{Proof of Theorem \ref{thm:privacy_guarantee} } \label{apx-thm:privacy_guarantee}

First, in the $k$th iteration of DP-PS, for all $z \in \mathcal{Z}$ and $i \in C(z)$, we denote by $P_{\mathcal{R}^k_{zi}(D_z)}(\tilde{y}^k_{zi})$ the probability density at any $\tilde{y}^k_{zi} \in S_k$, where $\mathcal{R}^k_{zi}$ is defined in \eqref{Randomized_fn_ACOPF} and $S_k$ is any subset of $\text{Range}(\mathcal{R}^k_{zi})$. Then we have

\vspace{-3mm}
\begin{small}
\begin{align*}
\mathbb{P} \{ \mathcal{R}^k_{zi} (D_z) \in S_k \}  = \int_{S_k} P_{\mathcal{R}^k_{zi}(D_z)}(\tilde{y}^k_{zi}) d\tilde{y}^k_{zi}.
\end{align*}
\end{small}
Now consider the following ratio:

\vspace{-3mm}
\begin{small}
\begin{align*}
  & \frac{P_{\mathcal{R}^k_{zi}(D'_z)}(\tilde{y}^k_{zi})}{P_{\mathcal{R}^k_{zi}(\bar{D}_z)}(\tilde{y}^k_{zi})}  = \frac{L( \tilde{y}^k_{zi} - \mathcal{Q}^k_{zi}(D'_z)| \bar{\Delta}^k_{zi} / \bar{\epsilon} )}{L(\tilde{y}^k_{zi} - \mathcal{Q}^k_{zi}(\bar{D}_z)| \bar{\Delta}^k_{zi}  / \bar{\epsilon})} \\
  & = \exp \Big( \big(\bar{\epsilon} / \bar{\Delta}^k_{zi}  \big) \big( \vert \tilde{y}^k_{zi} - \mathcal{Q}^k_{zi}(\bar{D}_z) \vert - \vert \tilde{y}^k_{zi} - \mathcal{Q}^k_{zi}(D'_z) \vert \big) \Big) \\
  & \leq \exp \Big( \big(\bar{\epsilon} / \bar{\Delta}^k_{zi}  \big) \big( \vert \mathcal{Q}^k_{zi}(D'_z) - \mathcal{Q}^k_{zi}(\bar{D}_z) \vert  \big) \Big) \\
  & \leq \exp(\bar{\epsilon}), \ \forall D'_z \in \widehat{\mathcal{D}}_{\beta}(\bar{D}_z),
\end{align*}
\end{small}
\noindent
where $L$ is from Definition \ref{def:laplace}, the first inequality holds due to the reverse triangle inequality, namely, $|a|-|b| \leq |a-b|$, and the last inequality holds since $\bar{\Delta}^k_{zi} \geq | \mathcal{Q}^k_{zi}(D'_z) - \mathcal{Q}^k_{zi}(\bar{D}_z) | $ for all $D'_z \in \widehat{\mathcal{D}}_{\beta}(\bar{D}_z)$ from \eqref{model:Sensitivity_problem}.
Similarly, one can obtain a lower bound as follows:

\vspace{-3mm}
\begin{small}
\begin{align*}  
  & \exp \Big( \big(\bar{\epsilon} / \bar{\Delta}^k_{zi}(\bar{\beta}) \big) \big( \vert \tilde{y}^k_{zi} - \mathcal{Q}^k_{zi}(\bar{D}_z) \vert - \vert \tilde{y}^k_{zi} - \mathcal{Q}^k_{zi}(D'_z) \vert \big) \Big) \\
  & \geq \exp \Big( -\big(\bar{\epsilon} / \bar{\Delta}^k_{zi}(\bar{\beta}) \big) \big( \vert \mathcal{Q}^k_{zi}(D'_z) - \mathcal{Q}^k_{zi}(\bar{D}_z) \vert  \big) \Big) \\
  & \geq \exp(-\bar{\epsilon}), \ \forall D'_z \in \widehat{\mathcal{D}}_{\beta}(\bar{D}_z),
\end{align*}
\end{small}
\noindent
where the first inequality holds due to the reverse triangle inequality, namely, $|a|-|b| \geq -|a-b|$.
Therefore, we have

\vspace{-3mm}
\begin{small}
\begin{align*}
  \exp(-\bar{\epsilon}) \leq \frac{P_{\mathcal{R}^k_{zi}(D'_z)}(\tilde{y}^k_{zi})}{P_{\mathcal{R}^k_{zi}(\bar{D}_z)}(\tilde{y}^k_{zi})} \leq \exp(\bar{\epsilon}), \ \forall D'_z \in \widehat{\mathcal{D}}_{\beta}(\bar{D}_z),
\end{align*}
\end{small}
\noindent
and integrating $\tilde{y}^k_{zi}$ over $S_k$ yields \eqref{epsilon_dp}.
This proves that $\bar{\epsilon}$-differential privacy on data is guaranteed for each iteration $k$ of DP-PS.

Second, for all $z \in \mathcal{Z}$ and $i \in C(z)$, we denote by $\mathcal{R}_{zi}$ a randomized function that maps the dataset $D_z \in \mathbb{R}^{|\mathcal{N}_z|}$ to $\tilde{y}_{zi} :=\{\tilde{y}^k_{zi} \}_{k=1}^K$, where $K$ is the total number of iterations consumed by DP-PS.
It suffices to show that

\vspace{-3mm}
\begin{small}
\begin{align}
\Big\vert \ln \Big( \frac{\mathbb{P} \{  \mathcal{R}_{zi} (D'_z) \in S \} }{ \mathbb{P} \{  \mathcal{R}_{zi} (\bar{D}_z) \in S \} } \Big) \Big\vert  \leq \bar{\epsilon}, \ \ \subalign{& \forall D'_z \in \widehat{\mathcal{D}}_{\beta}(\bar{D}_z), \\ & \forall S \subseteq \text{Range}(\mathcal{R}_{zi}).} \label{epsilon_diff_priv}
\end{align}
\end{small}
We denote by $P_{\mathcal{R}_{zi}(D_z)} (\tilde{y}_{zi})$ the joint density at any $\tilde{y}_{zi} \in S$.
Then we have

\vspace{-3mm}
\begin{small}
\begin{align*}
  \mathbb{P} \{  \mathcal{R}_{zi} (D_z) \in S \}
  & = \int_{S} P_{ \mathcal{R}_{zi}(D_z) } (\tilde{y}_{zi} ) d\tilde{y}_{zi}.
\end{align*}
\end{small}
The joint density function can be expressed by the conditional density functions:

\vspace{-3mm}
\begin{small}
\begin{align*}
& P_{ \mathcal{R}_{zi}(D_z) } (\tilde{y}_{zi} )=  P_{\mathcal{R}^1_{zi}(D_z), \ldots, \mathcal{R}^K_{zi}(D_z)} (\tilde{y}^1_{zi}, \ldots, \tilde{y}^K_{zi} ) \\
= & P_{\mathcal{R}^K_{zi}(D_z) |\mathcal{R}^1_{zi}(D_z), \ldots, \mathcal{R}^{K-1}_{zi}(D_z) } (\tilde{y}^K_{zi} | \tilde{y}^1_{zi}, \ldots , \tilde{y}^{K-1}_{zi}) \times \ldots \\
& \times P_{\mathcal{R}^1_{zi}(D_z) } (\tilde{y}^1_{zi}) \\
= &  \prod_{k=1}^K  L( \tilde{y}^k_{zi} - \mathcal{Q}^k_{zi}(D_z) | K \bar{\Delta}^k_{zi} / \bar{\epsilon} ) \\
= &  \prod_{k=1}^K \frac{1}{2K\bar{\Delta}^k_{zi}/\bar{\epsilon} } \exp(-\frac{|\tilde{y}^k_{zi} - \mathcal{Q}^k_{zi}(D_z) |}{K\bar{\Delta}^k_{zi}/\bar{\epsilon}} ).
\end{align*}
\end{small}
Taking similar steps in the first part of this proof, we obtain

\vspace{-3mm}
\begin{small}
\begin{align*}
\exp (-\bar{\epsilon}) \leq  \frac{P_{ \mathcal{R}_{zi}(D'_z) } (\tilde{y}_{zi} )}{P_{ \mathcal{R}_{zi}(\bar{D}_z) } (\tilde{y}_{zi} )} \leq \exp (\bar{\epsilon}), \ \forall D'_z \in \widehat{\mathcal{D}}_{\beta}(\bar{D}_z),
\end{align*}
\end{small}
and integrating $\tilde{y}_{zi}$ over $S$ yields \eqref{epsilon_diff_priv}.
This completes the proof.

\section{Proof of Theorem \ref{thm:Rule1_convergence_1}} \label{apx-thm:Rule1_convergence_1}
Since $H(\lambda^*) - \mathbb{E}[H_{\text{\tiny best}}(\lambda^k)] \geq 0$ and the right-hand side of \eqref{Rule1_sandwich_inequality} goes to $0$ as $K \rightarrow \infty$, \eqref{Rule1_convergence_expectation} holds.
Also, \eqref{Rule1_convergence_probability} holds due to Markov's inequality, namely, for $\epsilon > 0$,

\vspace{-3mm}
\begin{footnotesize}
\begin{align}
\mathbb{P} \big\{ {H}(\lambda^{\star}) - {H}_{\text{best}}(\lambda^{K}) \geq \epsilon \big\}  \leq \mathbb{E} \big[ {H}(\lambda^{\star}) - {H}_{\text{best}}(\lambda^{K}) \big] / \epsilon, \label{Markov}
\end{align}
\end{footnotesize}
\noindent
where the right-hand side of \eqref{Markov} goes to $0$ as $K \rightarrow \infty$.
From the right-hand side of \eqref{Rule1_sandwich_inequality}, the rate of convergence in expectation is $\mathcal{O}(G^{\text{\tiny U}}(\bar{\epsilon})/\log(K))$.
This completes the proof.

\section{Proof of Theorem \ref{thm:Rule1_convergence_as}} \label{apx-thm:Rule1_convergence_as}
By taking a conditional expectation on \eqref{Rule1_ineq}, we obtain

\vspace{-3mm}
\begin{small}
\begin{align*}
& \mathbb{E} \big[ \| \lambda^{k+1}-\lambda^{\star} \|^2 | \lambda^1, \ldots, \lambda^k \big]  \leq \| \lambda^k - \lambda^{\star} \|^2 + \alpha_k^2 G^{\text{\tiny U}}(\bar{\epsilon}),
\end{align*}
\end{small}
\noindent
where the inequality holds since $\alpha_k (H(\lambda^k) - {H}(\lambda^{\star}) ) \leq 0$ and $\mathbb{E}[\tilde{\xi}^k_{zi} | \lambda^k] = 0, \forall z \in \mathcal{Z}, \forall i \in C(z)$.
Since $\lambda^1$ is bounded by Assumption \ref{assump:sharp},
$\alpha_k^2 G^{\text{\tiny U}}(\bar{\epsilon}) \geq 0$, and $G^{\text{\tiny U}}(\bar{\epsilon}) \sum_{k=1}^{\infty} \alpha_k^2 < \infty$, the sequence $\{\lambda^k \}$ generated by Algorithm \ref{algo:DPPSA} with Rule 1 is a stochastic quasi-Feyer sequence for a set $\Lambda^{\star}$ of maximizers.
Based on Theorem 6.1 in \cite{ermoliev1988numerical} and the existence of a subsequence $\{\lambda^{k_s}\}$ such that ${H}_{\text{best}}(\lambda^{k_s})$ converges to ${H}(\lambda^{\star})$ with probability $1$ due to \eqref{Rule1_convergence_probability}, one can conclude that the sequence $\{\lambda^k \}$ converges to a point in $\Lambda^{\star}$. For more details, we refer the reader to the proof of Theorem 6.2 in \cite{ermoliev1988numerical}.

\section{Proof of Theorem \ref{thm:Rule2_convergence}} \label{apx-thm:Rule2_convergence}
Under Assumption \ref{assump:sharp}, it follows from \eqref{Rule2_basic_ineq_1} that

\vspace{-3mm}
\begin{small}
\begin{align}
  & \| \lambda^{k+1}-\lambda^{\star} \|^2  \leq \| \lambda^k - \lambda^{\star} \|^2 -  \frac{ \big( {H}(\lambda^{\star}) - {H}(\lambda^{k}) \big)^2}{\| \tilde{y}^k  \|^2}  + \nonumber \\
  &   2 \frac{{H}(\lambda^{\star}) - {H}(\lambda^{k})}{\| \tilde{y}^k  \|^2} \| s^k(\tilde{y}^k)- y^k \| \cdot \|\lambda^k-\lambda^{\star} \|  \nonumber \\
  & \leq \| \lambda^k - \lambda^{\star} \|^2 - \Big(1 -   \frac{2\| s^k(\tilde{y}^k)-{y}^k \| }{\mu} \Big) \frac{ \big( {H}(\lambda^{\star}) - {H}(\lambda^{k}) \big)^2}{\| \tilde{y}^k  \|^2} \nonumber \\
  & \leq \| \lambda^k - \lambda^{\star} \|^2 -  \big(G^{\text{\tiny L}} / G^{\text{\tiny U}}(\bar{\epsilon}) \big)  \big( {H}(\lambda^{\star}) - {H}(\lambda^{k}) \big)^2 ,
     \label{Rule2_basic_ineq_2}
\end{align}
\end{small}
\noindent
where the first inequality holds since $\tilde{\xi}^k = s^k(\tilde{y}^k)-{y}^k$, the second inequality holds due to Assumption \ref{assump:sharp}, and the last inequality holds due to $(1-2\|s^k-\hat{y}^k\|/\mu) \in (0,1]$ from Assumption \ref{assump:sharp} and Lemma \ref{lemma:Bound_y_tilde}.
By taking similar steps in Section \ref{sec:Rule2}, we obtain

\vspace{-3mm}
\begin{small}
\begin{align}
  0 \leq {H}(\lambda^{\star}) - \mathbb{E} \big[ \max_{k \in [K]} {H}(\lambda^{k}) \big]  \leq \sqrt{  \frac{ \lambda^{\text{\tiny U}} G^{\text{\tiny U}} (\bar{\epsilon}) }{  G^{\text{\tiny L}} K } }. \label{Rule2_sandwich_ineq_2}
\end{align}
\end{small}
Taking similar steps in the proof of Theorem \ref{thm:Rule1_convergence_1}, we conclude from \eqref{Rule2_sandwich_ineq_2} that the sequence produced by DP-PS with Rule 2 under Assumption \ref{assump:sharp} converges in expectation and in probability.
Also, the rate of convergence in expectation is $\mathcal{O}(G^{\text{\tiny U}}(\bar{\epsilon})/K)$.
It follows from \eqref{Rule2_basic_ineq_2} that  $\| \lambda^{k+1}-\lambda^{\star} \|^2  \leq \| \lambda^k - \lambda^{\star} \|^2$.
By taking a conditional expectation, we obtain

\vspace{-3mm}
\begin{small}
\begin{align*}
& \mathbb{E} \big[ \| \lambda^{k+1}-\lambda^{\star} \|^2 | \lambda^1, \ldots, \lambda^k \big]  \leq \| \lambda^k - \lambda^{\star} \|^2.
\end{align*}
\end{small}
Thus, the sequence $\{\lambda^k\}$ generated by DP-PS with Rule 2 under Assumption \ref{assump:sharp} is a stochastic quasi-Feyer sequence for a set $\Lambda^{\star}$ of maximizers.
As discussed in the proof of Theorem \ref{thm:Rule1_convergence_as}, it proves the convergence with probability $1$.
This completes the proof.

\section{DP-ADMM} \label{apx-DP-ADMM}
We present the augmented Lagrangian dual problem given by

\begin{scriptsize}
\begin{align*}
\max_{\lambda} \ \min \ & \sum_{z \in \mathcal{Z}} \Big\{ f_z(x_z) + \sum_{i \in C(z)} \Big( \lambda_{zi}( \phi_i - y_{zi}) + \frac{\rho}{2}( \phi_i - y_{zi})^2 \Big) \Big\}  \\
\mbox{s.t.} \ & (x_z, y_z) \in \mathcal{F}_z(\bar{D}_z), \ \forall z \in \mathcal{Z}, \\
& \phi_i \in \mathbb{R}, \ \forall i \in \mathcal{C},
\end{align*}
\end{scriptsize}
where $\lambda$ is a dual vector associated with constraint \eqref{model:ACOPF_matrix-2}.

For every iteration $k$ of the ADMM algorithm, it updates $(y^k, \phi^k, \lambda^k) \rightarrow (y^{k+1}, \phi^{k+1}, \lambda^{k+1})$ by solving a sequence of the following subproblems:

\begin{scriptsize}
\begin{subequations}
\begin{align}
& y^{k+1}_z \leftarrow  \argmin_{(x_z, y_z) \in \mathcal{F}_z(\bar{D}_z)} \ f_z(x_z) + \sum_{i \in C(z)} \Big\{ - \lambda^{k}_{zi} y_{zi} + \nonumber \\
& \hspace{30mm} \frac{\rho}{2}(\phi^{k}_{i}-y_{zi})^2\Big\}, \ \forall z \in \mathcal{Z}, \label{ADMM-1} \\
& \phi^{k+1}_i \leftarrow  \argmin_{\phi_i} \sum_{z \in F(i)} \lambda^{k}_{zi} \phi_i + \frac{\rho}{2}(\phi_{i} - y^{k+1}_{zi})^2, \forall i \in \mathcal{C}, \label{ADMM-2} \\
& \lambda^{k+1}_{zi} = \lambda^{k}_{zi} + \rho (\phi^{k+1}_{i} - y^{k+1}_{zi}), \ \forall z \in \mathcal{Z}, \forall i \in C(z).
\end{align}
\end{subequations}
\end{scriptsize}

In Algorithm \ref{algo:DPADMM}, we describe DP-ADMM.
In line 3, we solve the subproblem \eqref{ADMM-1} to find $y^{k+1}$, which is perturbed by adding the Laplacian noise $\tilde{\xi}^k$ in line 7.
In line 9, we solve the subproblem \eqref{ADMM-2} with $\tilde{y}^{k+1}$ to find $\phi^{k+1}$.
In line 10, the dual variable is updated from $\lambda^k$ to $\lambda^{k+1}$.

\begin{algorithm}[H] 
  \caption{DP-ADMM}
  \begin{algorithmic}[1]
  \STATE Set $k \leftarrow 1$, $\lambda^1 \leftarrow 0$, and $\phi^1 \leftarrow 0$.
  \FOR{$k \in \{1, \ldots, K\}$}
  \STATE
  Given $\lambda^k$ and $\phi^k$, find $y^{k+1}$ by solving \eqref{ADMM-1}.  
  \STATE \textbf{\# Perturbation of $y^{k+1}$}
  \STATE Solve \eqref{model:Sensitivity_problem} to find $\{\bar{\Delta}^k_{zi} \}_{z \in \mathcal{Z}, i \in C(z)}$.
  \STATE Extract $\tilde{\xi}^k_{zi}$ from
  $L(\tilde{\xi}^k_{zi} \vert \bar{\Delta}^k_{zi}  / \bar{\epsilon})$ in Definition \ref{def:laplace}.
  \STATE $\tilde{y}^{k+1}_{zi} \leftarrow y^{k+1}_{zi}+\tilde{\xi}^k_{zi}, \ \forall z \in \mathcal{Z}, \forall i \in C(z).$.
  \STATE \textbf{\# Solve the second-block problem}  
  \STATE Solve \eqref{ADMM-2} with $\tilde{y}^{k+1}$.
  \STATE \textbf{\# Update dual variables}  

  \vspace{-3mm}
  \begin{small}
  \begin{align*}
    \lambda^{k+1}_{zi} = \lambda^{k}_{zi} + \rho (\phi^{k+1}_{i} - \tilde{y}^{k+1}_{zi}), \ \forall z \in \mathcal{Z}, \forall i \in C(z). 
  \end{align*} 
  \end{small} 
  \ENDFOR
  \end{algorithmic}
  \label{algo:DPADMM}
\end{algorithm}



\end{document}